\RequirePackage{fix-cm}
\documentclass[11pt,a4paper]{article}

\usepackage{bm}
 
\usepackage{xspace}
\usepackage{geometry}
\usepackage{array}
\usepackage{subfigure}
\usepackage{color}
\usepackage{epsfig}
\usepackage{hhline}
\usepackage{natbib}
\usepackage{bbm}
\usepackage[T1]{fontenc}
\usepackage[]{graphicx,float,latexsym,times}
\usepackage{amsfonts,amstext,amsmath,amssymb}
\usepackage[utf8]{inputenc}    
\DeclareMathAccent{\widehat}{\mathord}{largesymbols}{"62}
\DeclareMathAccent{\widetilde}{\mathord}{largesymbols}{"65}

\def\sign{{\mathrm{sgn}}}

\def\ebo{\textrm{\mathversion{bold}$\mathbf{\beta}^0$\mathversion{normal}}}

\def\oo{\textrm{\mathversion{bold}$\mathbf{0}$\mathversion{normal}}}
\def\eb{\textrm{\mathversion{bold}$\mathbf{\beta}$\mathversion{normal}}}

\def\eo{\textrm{\mathversion{bold}$\mathbf{\omega}$\mathversion{normal}}} 
\def\ek{\textrm{\mathversion{bold}$\mathbf{\kappa}$\mathversion{normal}}}

\def\eXL{\textrm{\mathversion{bold}$\mathbf{\Xi}$\mathversion{normal}}} 

\def\eE{\mathbb{E}}

\def\e1{1\!\!1}

\newcommand{\SSS}{\mathbb{S}}

\newtheorem{theorem}{Theorem}

\newtheorem{remark1}{Remark} 

\def\ee1{\textrm{\mathversion{bold}$\mathbf{\varepsilon}$\mathversion{normal}}}

\def\oo{\textrm{\mathversion{bold}$\mathbf{0}$\mathversion{normal}}}

\def\eu{\mathbf{{u}}}

\newcommand{\N}{\mathbb{N}}
\newcommand{\R}{\mathbb{R}}
\newcommand{\PP}{\mathbb{P}}

\def\eX{\mathbf{X}}

\newcommand{\Var}{\mathbb{V}\mbox{ar}\,}

\def\argmin{\mathop{\mathrm{arg\,min}}} 

\begin{document}
 
\title {{Right-censored models by the expectile method}}
 
\date{}
\maketitle
\author{
	\begin{center}
		Gabriela CIUPERCA \\
		\small{Université Claude Bernard Lyon 1, UMR 5208, Institut Camille Jordan, \\
			Bat.  Braconnier, 43, blvd du 11 novembre 1918, F - 69622 Villeurbanne Cedex, France.
		}
	\end{center}
}
  
\begin{abstract}
Based on the expectile loss function and the adaptive LASSO penalty, the paper proposes and studies the estimation methods for the accelerated failure time (AFT) model. In this approach, we need to estimate the survival function of the   censoring variable   by  the Kaplan-Meier estimator.  The  AFT model parameters are first estimated by the expectile method and afterwards, when the number of explanatory variables can be large, by the adaptive LASSO expectile method which directly carries out the automatic selection of variables. We also obtain the convergence rate and asymptotic normality for the two estimators, while showing the sparsity property for the censored adaptive LASSO expectile estimator. A numerical study using Monte Carlo simulations confirms the theoretical results and demonstrates the competitive performance of the two proposed estimators. The usefulness of these estimators is illustrated by applying them to three survival data sets. \\

 \noindent \textbf{Keywords}: {Accelerated failure time  \and right-censoring \and expectile \and LASSO \and automatic selection \and asymptotic behavior.}\\
 \noindent \textbf{MSC 2020 Subject Classification}: {62N01, 62N02, 62F12, 62J07.}
\end{abstract}

\section{Introduction}
Censored models have many practical applications, for example in medicine, industry and economics, to name but a few of these areas. As highlighted in \cite{Wei.92}, in survival analysis, the accelerated failure time (AFT) model is a useful alternative  to the Cox model. This usefulness stems from the fact that the survival time is the response variable in a linear regression model, which  encourages the consideration for the AFT model of the estimation techniques originally used for  linear regression models.  For the properties and classical estimation methods of an AFT model, interested readers can refer to the book by \cite{Kalbfleisch.Prentice.02}. Various estimation techniques for these models have been considered in the very rich literature on the subject. 
Historically, the parameters of the linear model associated with a censored model have been estimated using the least squares (LS) technique (see \cite{Ritov.90}, \cite{Tsiatis.90}, \cite{Li.Wang.12}, \cite{Jin.Jung.Wei.06} among others). If the model errors do not satisfy the classical conditions, then the LS estimation is sensitive to outliers. One possible approach in this case was to consider censored median regression models (see \cite{Zhou.Wang.05}, \cite{Zhou.06}, \cite{Huang.Ma.Xie.07}), which were afterwards generalized by quantile models. 
In this context,  \cite{Portnoy.03}, \cite{Wang.Wang.09} propose a censored    weighted quantile estimator, in the latter paper, consistency and asymptotic normality are shown. The same asymptotic properties are satisfied by the estimator proposed in \cite{Peng.Huang.08} by a new quantile regression approach for survival data subject to conditional independent censoring. A censored quantile regression with the explanatory variables  measured with errors is studied by \cite{Ma.Yin.11} who propose a composite objective function based on inverse censoring-probability weighting. The  obtained estimator improves the efficiency of the estimation. The same model is studied by \cite{Wu.Ma.Yin.15} who propose a smoothed martingale estimating equation and,   for the estimation of the model parameters, generalize the grid-based estimation procedure proposed by \cite{Peng.Huang.08}. \cite{De Backer.19} considered an alternative approach, by adapting the loss function with an estimator of the survival function of the censoring variable. They proposed an algorithm to minimize the adapted loss function, resulting in a consistent and asymptotically normal estimator.  \\
Note that, in survival models, it is necessary to use a consistent estimator to estimate the distribution of the censoring variable. The most popular estimator is the Kaplan-Meier estimator, whose properties can be found in \cite{Stute.94} or  \cite{Wang.Ng.08}.
In practical problems, especially when the model has a large number of explanatory variables, it is often necessary to automatically select the relevant variables. This can be achieved by applying an adaptive LASSO penalty to the loss function. This penalty was originally introduced and extensively studied for linear models with uncensored response variables. For more information, please refer to the papers by  \cite{Zou.06}, \cite{Wu:Liu:09}, \cite{Xu:Ying:10}, \cite{Liao-Park-Choi.19} and \cite{Ciuperca.21}. In the context of survival data analysis, the adaptive LASSO penalty is considered by \cite{Johnson.09},  \cite{He.Wang.Zhou.19} for a censored LS model, by \cite{Shows.Lu.Zhang.10} for a censored median model, by  \cite{Tang-Zhou-Wu.12}, \cite{Zheng.Peng.He.18}, \cite{Wang.Jiang.Xu.Fu.2021} for a censored quantile model. In a high-dimensional censored model,  \cite{Lee.Park.Lee.23} used the quantile forward method with an extended BIC penalty to select significant variables.
The quantile estimation method has the disadvantage that the loss function is not differentiable. This is a problem in theoretical studies and for related numerical methods. A possible solution is  to consider the expectile estimation method, introduced by \cite{Newey-Powell.87} for a classical linear model,  as a generalization of the LS method. To the best of the author's knowledge, the expectile estimation method has received little attention in the literature on censored models. However,  \cite{Seipp.Uslar.21} considers it under the assumption that the distribution function of the censoring variable is known and only states the asymptotic normality of the corresponding censored expectile estimator. Furthermore,  the topic of automatic variable selection is not addressed.  Not last, also of relevance for the present paper, when the censored model variables are grouped, to remove the unimportant groups, \cite{Huang.Kopciuk.Lu.20} consider the censored adaptive group bridge LS method. However,  \cite{Li.Gu.12}, used the penalized log-marginal likelihood method to perform variable selection for general transformation models with right-censored data, using the adaptive LASSO penalty. Likewise,  \cite{Cai.Huang.09} proposed a rank-based adaptive LASSO estimator for a AFT model with high-dimensional predictors, while  \cite{Chung.Long.Johnson.13} reviewed regularized rank-based coefficient estimation procedures.  Other penalties for variable selection in a censored models have also been considered in the literature. \cite{He.Zhou.Zou.20} investigated a sparse and consistent estimator for length-based data using the SCAD penalty, while \cite{Huang.Ma.10} modeled the relationship between covariates and survival  time using the accelerated failure time models, with bridge penalization. The  weighted LS loss function penalized with a minimax concave penalty induces an estimator which enjoys oracle property in \cite{Hu.Chai.13}. \cite{Wang.Song.2011} showed that the adaptive LASSO of the weighed LS estimator in the AFT model with multiple covariates has oracle properties and used the BIC criterion for tuning parameter selection.  Always a weighted LS loss function combined with $L_0$-penalization or a seamless $L_0$-penalization is considered in \cite{Cheng.Feng.Huang22}, \cite{Xu.Wang.23}, respectively. Notable works include those of \cite{Stute.93}, \cite{Stute.96} who consider the Kaplan-Meier weights for the censored weighted LS estimator. This estimator is   $n^{-1/2}$-consistent and asymptotically normal. This estimator was subsequently used by \cite{Su.Yin.Zhang.Zhao.23} as an adaptive weight to define the censored adaptive LASSO LS estimator  for an AFT model with an exceptionally large sample size, where the dimension of the explanatory variables is large but smaller than the sample size. \cite{Su.Yin.Zhang.Zhao.23} develop a divide-and-conquer approach based on the censored adaptive LASSO LS estimator, which will generate an estimator for which the oracle properties are proven.   \\
In this paper, we consider a random right-censoring  model estimated by the expectile method. We also propose an estimator that allows the automatic selection of significant explanatory variables.   We emphasize that the numerical study and applications on real data presented in this paper show the practical interest of the two estimation methods for censored models with respect to other methods presented in the literature. More specifically, the main contributions of this paper are fourfold. First, we introduce the censored expectile estimator for an AFT model when the errors may be asymmetric. Theoretically, we prove that this estimator is $n^{-1/2}$-consistent and asymptotically normal. Second, we introduce and study the censored adaptive LASSO expectile estimator which is interesting and useful for automatic variable selection  in an AFT model with a large number of explanatory variables. Theoretically, we prove the $n^{-1/2}$-consistency of this estimator and its oracle properties.  Third, we confirm the theoretical results of the two estimators and demonstrate their competitive performance through a numerical study. Fourth, we illustrate the usefulness of these two censored estimators on three survival applications.  \\
The paper is organized as follows.  The model, assumptions and general notations are introduced in  Section \ref{section_model}. In  Section \ref{section_results} we introduce  the censored expectile estimator and study its theoretical properties. Afterwards, the censored adaptive LASSO expectile estimator is defined, followed by the study of its asymptotic behavior, whose oracle properties. Section \ref{section_simus} presents the simulation results, followed by three applications on survival datasets in Section \ref{section_application}. The proofs of the theoretical results are relegated in Section \ref{section_proofs}.
  
\section{Model and assumptions}
\label{section_model}
In this section we present the statistical model, the necessary assumptions and we also introduce some of the notations used throughout in the paper.\\
Let's start with some notation that will be used throughout the rest of the paper. Note that all vectors are considered column. Moreover, matrices and vectors are denoted by boldface uppercase and lowercase  letters, respectively. For a vector or matrix, the symbol $T$ at the top right is used for its transpose. We denote the Euclidean norm of a vector  by  $\|.\|$  and the $q$ vector with all components 0 by $\textbf{0}_q$. For an event $E$, $\e1_E$ denotes the indicator function that the event $E$ occurs.  Given a set ${\cal S}$, we denote its cardinality by $|{\cal S}|$ and its complementary set by ${\cal S}^c$.   For a real $x$ we use the notation $\textrm{sgn}(x)$ for the sign function $\textrm{sgn}(x)={x}/{|x|}$: if $x \neq 0$ and $\textrm{sgn}(0)=0$. We use $ \overset{\PP} {\underset{n \rightarrow \infty}{\longrightarrow}} $, $ \overset{\cal L} {\underset{n \rightarrow \infty}{\longrightarrow}} $ to represent the convergence in probability and in distribution, respectively, as $n$ converges to infinity.  The following notations are also used throughout in the paper: if $U_n$ and $V_n$ are two random variable sequences,  the notation $V_n=o_\PP(U_n)$ means that $\underset{n \rightarrow \infty}{\text{lim}} \mathbb{P}(|  {U_n}/{V_n} | >e)=0 $ for all $ e >0$ , while  the notation $V_n=O_\PP(U_n)$  means that there exists $ c >0 \, \text{so that} \, \underset{n \rightarrow \infty}{\text{lim}} \mathbb{P}(|  {U_n}/{V_n} | >c )<e  $ for all $ e>0$. We also note by ${\cal N}_p(\textbf{m},\textbf{V})$ a multivariate normal distribution of dimension $p$, with the mean a $p$-vector $\textbf{m}$ and  the  variance matrix $\textbf{V}$.\\

Consider the following linear model on  $n$ observations:
\begin{equation}
\label{eq1}
T^*_i=\eX_i^\top \eb+\varepsilon_i, \qquad i=1, \cdots , n, 
\end{equation}
where $\eX_i$ is a  random vector of $p$ observable explanatory variables,  $\varepsilon_i$ is a  continuous random variable of errors,  $\eb \in \R^p$ is the vector of parameters and $\ebo=(\beta^0_1, \cdots , \beta^0_p)$ its true  (unknown) value.  Taking into account the unobservable  random variable  $T^*_i$, we consider  $T_i=\exp(T_i^*)$, which is the failure time (or survival time).\\
Let   ${\cal C}_i$ be the censoring variable (censoring time) for the $i$th observation.  In this paper we assume that $T_i$ is randomly right-censored by ${\cal C}_i$, such that the censoring time $T_i$ cannot always be observed. Thus, due to censoring, the observed variables are $(Y_i,\eX_i,\delta_i)_{1 \leqslant i \leqslant n}$, with $Y_i \equiv \min(T_i,{\cal C}_i)$ and the corresponding failure indicator $\delta_i \equiv \e1_{T_i \leq {\cal C}_i}$.  Note that, the random variable $\delta_i$ indicates whether $T_i$ was observed or not, being the censoring indicator. Recall that this type of model is called an accelerated failure time (AFT) model.\\
The focus of this paper is on inference for the parameter $\eb$ of model (\ref{eq1}).\\

In the following, throughout this paper we will denote by $c$ a generic constant, without interest, which does not depend on $n$. We also denote by $Y$, ${\cal C}$, $T$, $\eX$, the generic variable for $Y_i$, ${\cal C}_i$, $T_i$, $\eX_i$. The components of  $\eX$ are $(X_1, \cdots , X_p)$ and those of  $\eX_i$ are $(X_{1i}, \cdots , X_{pi})$.\\
For all  $t>0$ we define $G_0(t) \equiv \PP_{\cal C}[{\cal C} >t]$, which is the survival function of the censoring variable ${\cal C}$. Then, $1-G_0(t)$ is the distribution function of ${\cal C}$. \\
We denote by $\PP_\eX$ the probability law of the random vector $\eX$ and by $\eE_\eX$ the expectation with respect to the distribution of $\eX$.  Similarly, we denote by $\PP_{\cal C}$, $\PP_\varepsilon$, the probability laws of $\cal C$ and $\varepsilon$, respectively and by $\PP$ the joint probability of $(\eX, {\cal C}, \varepsilon)$, and $\eE$ the corresponding expectation. \\

Often in applications, the survival function $G_0$  of the censoring variable $\cal C$ is unknown. There are several estimators for $G_0$, including the best known: Kaplan-Meier and Fleming-Harrington estimators (see  \cite{Stute.94}, \cite{Wang.Ng.08} for the Kaplan-Meier estimator proprieties and \cite{Fleming.Harrington.84} for those of the Fleming-Harrington estimator). In this paper we consider the Kaplan-Meier estimator. For $t >0$,   the Kaplan-Meier estimator of  $G_0(t)$ based on the random variables $(Y_i,\delta_i)_{1 \leqslant i \leqslant n}$ is defined by:
\[
\widehat G_n(t) \equiv \prod^n_{ \substack{
		i=1\\ Y_i \leq t}} \bigg( \frac{n -R_i}{n-R_i+1}\bigg)^{\delta_i},
\]
with $R_i$ the rank of $Y_i$ in $(Y_i)_{1 \leqslant i \leqslant n}$.\\

We now present the classical assumptions  on  ${\cal C}_i$, $\varepsilon_i$, $\eX_i$, $T_i$, for any $i=1, \cdots , n$. They will be needed in this paper.
\begin{description}
	\item \textbf{(A1)}  ${\cal C}_i$ is independent of  $\varepsilon_i$  and of $\eX_i$.  The  censoring variable ${\cal C}_i$ is also independent of the  failure time $T_i$  conditional on $\eX_i$.
\item \textbf{(A2)} The random vectors $(T_i,{\cal C}_i,\eX_i)_{ 1 \leqslant i \leqslant n}$ are independent and identically distributed (i.i.d.).
\item  \textbf{(A3)} The random variables $(\varepsilon_i)_{ 1 \leqslant i \leqslant n}$ are i.i.d.
\item  \textbf{(A4)} $\PP[t \leq T \leq {\cal C}] \geq \zeta_0 >0$ for all $t \in [0,B]$ with  $\zeta_0$ a positive  constant and  constant $B$  the maximum follow-up.
\item  \textbf{(A5)} The random vector $\eX$ is bounded (there exists $c>0$ such that $\PP_\eX [\| \eX\| < c]=1$) and $n^{-1} \sum^n_{i=1} \eX_i \eX_i^\top \overset{\PP_\eX} {\underset{n \rightarrow \infty}{\longrightarrow}}  \eE_\eX [\eX \eX^\top]$, with $ \eE_\eX [\eX \eX^\top]$ a  positive definite matrix.
\end{description}

\medskip
For the survival function $G_0$ we consider the classical condition:
\begin{description}
	\item \textbf{(A6)}  $G_0$ is continuous and its derivative is uniformly bounded on $[0,B]$.
\end{description}

\medskip
Moreover, we obviously have $\eE_{\cal C}[\delta_i]=\eE_{\cal C}[\e1_{T_i \leq {\cal C}_i}]=\PP_{\cal C}[{\cal C}_i \geq T_i]=G_0(Y_i)$ for any $i=1, \cdots ,n$.\\  
These assumptions are commonly considered in the literature. Note that  \cite{Wang.Jiang.Xu.Fu.2021}  considers (A2), (A3) and that  $T_i$ and ${\cal C}_i$  are independent conditional on $\eX_i$. Moreover,   \cite{Tang-Zhou-Wu.12}  also assumes (A1), while  \cite{Zhou.06} considers ${\cal C}_i$  independent of $T_i$  conditioned by $\eX_i$, together with the facts that $\varepsilon_i$ are i.i.d., the survival function  $G_0$ does not depend on  $\eX_i$ and assumption (A3). Assumption (A3) is also considered by \cite{Li.Wang.12}  which also assumes that $Y_i$ is independent of ${\cal C}_i$ conditioned by $\eX_i$. \cite{Seipp.Uslar.21}, which also considers a right-censored model by the expectile method, supposes that ${\cal C}_i$ is independent of $T_i$ conditioned on $\eX_i$ and that ${\cal C}_i$ is independent of $\eX_i$. The paper by  \cite{Shows.Lu.Zhang.10} assumes that ${\cal C}_i$  is independent of $T_i$ and of $\eX_i$.  \cite{Sun.Zhang.09} and \cite{De Backer.19} take assumption (A2) , ${\cal C}_i$ independent of  $T_i$ conditioned by  $\eX_i$. \cite{Portnoy.03} considers, in addition to (A2), that the distribution of ${\cal C}_i$ may depend on $\eX_i$ but that, conditioned by $\eX_i$, the random variables $T_i$ and ${\cal C}_i$ are independent. \cite{Wang.Wang.09} assumes that the distribution of ${\cal C}_i$ depends on $\eX_i$,  together with  assumption (A2).  \cite{Zhou.Wang.05} assumes (A3), together with ${\cal C}_i$ independent of $\varepsilon_i$. \cite{Ying.Jung.Wei.95}  assume that $T_i$ and ${\cal C}_i$ are independent, (A2), (A3). Assumption (A4) was also considered by \cite{Lee.Park.Lee.23}, \cite{Tang-Zhou-Wu.12} and \cite{Shows.Lu.Zhang.10} for the consistency of the Kaplan-Meier estimator $\widehat G_n$. Assumption (A5) is also considered by \cite{Zhou.06} for censored median regression, while (A6) is used in the paper of \cite{Ying.Jung.Wei.95}.\\
Note that assumptions (A4) and (A6) are necessary for the consistency of $\widehat G_n$ and for the Taylor expansion of $\widehat G_n(t)$. \\
 
To estimate the parameter vector $\eb$ in the basis of $(Y_i,\eX_i,\delta_i)_{ 1 \leqslant i \leqslant n}$ we consider the expectile function:
\[
\rho_\tau(x)=|\tau - \e1_{x <0}| x^2, \qquad x \in \R,
\]  
with $\tau \in (0,1)$ the expectile index.\\
The derivative of $\rho_\tau(x)$   is $g_\tau(x) \equiv\rho'_\tau(x-t)|_{t=0}= - 2 \tau x \e1_{x \geq 0}- 2(1-\tau)x \e1_{x<0}$ and the second derivative is $h_\tau(x) \equiv \rho''_\tau(x-t)|_{t=0}=2 \tau \e1_{x \geq 0}+2(1-\tau) \e1_{x<0}$. \\
The interest of the expectile estimation method is multiple. First of all it can be applied when the distribution of the model errors $\varepsilon$  is asymmetric, in which case  the LS estimation method is not accurate because the corresponding estimators are less efficient. One option is the quantile method but this has the disadvantage that the loss function is not derivable which complicates the theoretical study and the computational methods,  especially in the case of a censored model.\\

In order to study the properties of the estimators proposed  in this paper, we consider the following assumption for the errors $(\varepsilon_i)_{1 \leqslant i \leqslant n}$ in addition to (A3):
\begin{description}
	\item \textbf{(A7)}   $\eE_\varepsilon[\varepsilon^4] < \infty$ and $\eE_\varepsilon[g_\tau(\varepsilon)]=0$.
\end{description}

If $\tau=1/2$ then $\eE_\varepsilon[g_\tau(\varepsilon)]=0$ in assumption (A7) becomes $\eE_\varepsilon[\varepsilon]=0$, which is the standard condition considered for the LS estimation method.\\
Moreover, assumption (A7) is often required for the expectile models (see e.g. \cite{Gu-Zou.16}, \cite{Ciuperca.21}). 
Note that the second assumption of (A7) implies that the expectile index $\tau $ is fixed.\\

Before proceeding to the main theoretical results presented in  the following section, we introduce some random processes and random vectors.\\
For $s \in [0,B]$, let be the following random $p$-vector 
$$\ek(s) \equiv \lim_{n \rightarrow \infty} \frac{1}{n} \sum^n_{i=1}  \frac{\delta_i}{G_0(Y_i)}  \e1_{Y_i \geq s} \eX_i g_\tau(\varepsilon_i),$$
which is bounded on  $[0,B]$ by assumption (A4),  the constant $B$ being defined in assumption  (A4). For $s, t \in [0,B]$, $j=1, \cdots , n$, consider the following random processes:
\begin{equation*}
\left\{  
\begin{split}
y(s)  &\equiv \lim_{n \rightarrow \infty} \frac{1}{n} \sum^n_{i=1} \e1_{Y_i \geq s}, \\
M_j^{\cal C}(t) & \equiv (1- \delta_j) \e1_{Y_j \leq t} - \int^t_0 \e1_{Y_j \geq s} d \Lambda_{\cal C}(s),
\end{split}
\right. 
\end{equation*}
where $\Lambda_{\cal C}$ is the cumulative hazard function of the censoring  variable ${\cal C}$, i.e.  $\Lambda_{\cal C}(t) \equiv -\log(G_0(t))$ (see e.g. \cite{Sun.Zhang.09}, \cite{Shows.Lu.Zhang.10}, \cite{Tang-Zhou-Wu.12}, \cite{Wang.Jiang.Xu.Fu.2021}). \\
Let us remark that $\{ M_j^{\cal C}(t)\}$ is a  martingale with respect to the $\sigma$-filtration: $\sigma \{\e1_{Y_j \geq u}, (1-\delta_j) \e1_{Y_j \leq u}, 0 \leq u \leq t; \eX_j; j=1, \cdots, n\}$ (see \cite{Sun.Zhang.09}).
\section{Estimators and asymptotic behavior}
\label{section_results}
In this section we present our theoretical results. More specifically, we introduce and study two estimators for the parameter $\ebo$. First, we define the censored expectile estimator, find its convergence rate, and show its asymptotic normality.  Afterwards, we define an adaptive LASSO type estimator for which the asymptotic properties are studied. The proofs of the theorems presented in this section are relegated in Section \ref{section_proofs}.\\

Based on the consideration that the random variable $T_i$ is unobservable, together with the fact that the survival function $G_0$ is unknown, and  since model (\ref{eq1}) is right-censored, then we will consider as an estimator for  $\ebo$:
\[
\widetilde \eb_n \equiv \argmin_{\eb \in \R^p} \sum^n_{i=1} \frac{\delta_i}{\widehat G_n(Y_i)} \rho_\tau(\log(Y_i)- \eX_i^\top \eb ),
\]
that which we call the censored expectile estimator. The components of the $p$-vector $\widetilde \eb_n $ are $\big(\widetilde \beta_{n,1}, \cdots, \widetilde \beta_{n,p}\big)$ and $\widehat G_n$ is the Kaplan-Meier estimator of $G_0$. For the particular case $\tau=1/2$, the estimator $\widetilde \eb_n $ becomes the censored LS estimator. \\
A similar estimator was proposed by \cite{Ying.Jung.Wei.95}, \cite{Zhou.06}, who considered the  $L_1$ norm  instead of the expectile function $\rho_\tau$ and they get the censored median estimator. Later their results were  generalized by \cite{Wang.Wang.09}, \cite{Wang.Jiang.Xu.Fu.2021} who studied  the censored quantile estimators.\\

Our first theoretical result concerns the asymptotic behavior of the estimator $\widetilde{\eb}_n$. The following theorem  gives the convergence rate of the censored expectile estimator  and its asymptotic normality.
\begin{theorem}
\label{Theorem2.1 Tang}
Under assumptions (A1)-(A7) we have:\\
(i) $\widetilde \eb_n -\ebo=O_\PP(n^{-1/2})$.\\
(ii) $n^{1/2}(\widetilde \eb_n - \ebo)   \overset{\cal L} {\underset{n \rightarrow \infty}{\longrightarrow}} {\cal N}_p(\oo_p,\textbf{S}^{-1}_3(\textbf{S}_1+\textbf{S}_2)\textbf{S}^{-1}_3)$, 
with the $p$-square matrices: \\
$\textbf{S}_1\equiv \eE_\varepsilon\big[ g^2_\tau(\varepsilon)\big] \eE_\eX\big[ \eX \eX^\top/ G_0(Y)\big]$, $\textbf{S}_2 \equiv \eE \big[\int^B_0 {\ek(s) \ek^\top(s)}/{y(s)}d \Lambda_C(s) \big]$ and $\textbf{S}_3 \equiv  \eE_\varepsilon \big[h_\tau(\varepsilon)\big] \eE_\eX [\eX \eX^\top]$. 
\end{theorem}

\begin{remark1}
	Note that we have $\eE_\varepsilon[h_\tau(\varepsilon)] >0$ for all $\tau \in (0,1)$ and for $\tau=1/2$ (i.e. the LS method) we get  $\eE_\varepsilon[ g_\tau(\varepsilon)]=\eE_\varepsilon[  \varepsilon]=0$, $\eE_\varepsilon\big[ g^2_\tau(\varepsilon)\big]=\Var[\varepsilon]$, $\eE_\varepsilon[h_{0.5}(\varepsilon)] =1$. 
	\end{remark1}

The results of Theorem \ref{Theorem2.1 Tang} are a generalization of those obtained by \cite{Jin.Jung.Wei.06}, \cite{Johnson.09}, \cite{Li.Wang.12} for the censored LS estimator. Similar consistency and asymptotic normality results, but with different asymptotic variance-covariance matrices, have previously been obtained by \cite{Zhou.06}, \cite{Wang.Wang.09}, \cite{Ying.Jung.Wei.95}, \cite{Portnoy.03} for censored median or quantile models.\\
The interest of the result of Theorem \ref{Theorem2.1 Tang}\textit{(ii)} is to be able to construct the confidence interval for $\eb$ or to carry out hypothesis tests on the components of $\eb$.\\

The censored expectile estimator $\widetilde \eb_n$ allows us to build a new estimator that has the property of automatically selecting the relevant variables. This property is particularly useful when the number $p$ of explanatory variables is large. In the present paper we consider  $p< n$. Moreover, the number $p$ can be close to $n$, but it does not depend on $n$. Then, for censored model  (\ref{eq1}), we define  the censored adaptive  LASSO expectile estimator, as follows:
\begin{equation}
	\label{ebn}
\widehat \eb_n \equiv \argmin_{\eb \in \R^p} \bigg(  \sum^n_{i=1}\frac{\delta_i}{\widehat G_n(Y_i)}\rho_\tau(\log(Y_i)-\eX_i^\top \eb)+ \lambda_n \sum^p_{j=1} \widehat \omega_{n,j}|\beta_j| \bigg),
\end{equation}
with the adaptive weights $ \widehat \omega_{n,j}\equiv |\widetilde \beta_{n,j}|^{-\gamma} $ and $\gamma >0$ a known parameter. The estimator $\widehat \eb_n $ is written as $\big(\widehat \beta_{n,1}, \cdots , \widehat \beta_{n,p}\big)$.  The tuning parameter $(\lambda_n)_{n \in \N}$ is a positive deterministic sequence which together with the weights $ \widehat \omega_{n,j}$, controls the overall complexity of the model.  \\
Furthermore, we emphasize that for the particular case $\tau=1/2$, the estimator $\widehat \eb_n$ becomes the censored adaptive LASSO LS estimator, an estimator studied by \cite{Johnson.09} for $\gamma=1$.\\
In order to show the asymptotic properties of this estimator, let us introduce the index set of true non-zero coefficients of model (\ref{eq1}): 
$$
{\cal A} \equiv \{j \in \{ 1, \cdots ,p\}; \; \beta^0_j \neq 0\}.
$$
 Since $\ebo$ is unknown,, so is the set ${\cal A}$. Without reducing the generality, we assume that  ${\cal A}$ contains  the first  $q \equiv |{\cal A}|$ natural numbers: ${\cal A}=\{1, \cdots , q\}$.  So its complementary set is ${\cal A}^c=\{q+1, \cdots , p\}$. The adaptive penalty of optimization problem (\ref{ebn}) allows the sparse estimation of the coefficients. \\
 
By the following theorem we prove that the censored adaptive LASSO expectile estimator $\widehat \eb_n$ has the same convergence rate as the censored expectile estimator $\widetilde \eb_n$. This result shows that the penalty has no effect on the convergence rate.
\begin{theorem}
	\label{Theorem2.2 Tang}
Under assumptions (A1)-(A7), if the tuning parameter sequence $(\lambda_n)_{n \in \N}$ is such that  $n^{-1/2} \lambda_n =O_\PP(1)$, then $\widehat\eb_n -\ebo=O_\PP(n^{-1/2})$.
	\end{theorem}
The result of Theorem \ref{Theorem2.2 Tang} will be useful to show that the estimator $\widehat\eb_n$ satisfies the oracle property,  i.e. that it is sparse and that the estimators of the non-zero coefficients are asymptotically normal. To do this, similar to ${\cal A}$, we consider the index set of  non-zero estimated coefficients:  
$$
\widehat {\cal A}_n \equiv \{j \in \{ 1, \cdots , p\}; \; \widehat \beta_{n,j} \neq 0\},
$$
 which is an estimator of the set ${\cal A}$.\\
We denote by  $\eX_{{\cal A}}$ the sub-vector of $\eX$ with indices in  ${\cal A}$ and similarly $\ek_{\cal A}$ the sub-vector  of $\ek$.   We also denote the following sub-vectors of $\eX_i$ : $\eX_{{\cal A},i} \equiv (X_{ji})_{j \in {\cal A}}$ and  $\eX_{{\cal A}^c,i} \equiv (X_{ji})_{j \in {\cal A}^c}$. For a vector $\eb$, we use the notational conventions $\eb_{\cal A}$ for its sub-sector containing the corresponding components of ${\cal A}$. We also use  the notation for the following $|{\cal A}|$-vector:  $|\eb^0_{\cal A}|^{-\gamma}\equiv\big(|\beta^0_1|^{-\gamma}, \cdots, |\beta^0_q|^{-\gamma}\big)$. \\ 
Let be the random $|{\cal A}|$-vector 
$$\textbf{v}_3 \equiv \frac{\delta}{  G_0(Y)} g_\tau(\varepsilon)  \eX_{{\cal A}}+\int^B_0 \frac{\ek_{\cal A}(s)}{y(s)}dM^{\cal C}(s) 
$$
 and the $|{\cal A}|$-square matrix 
 $$
 \textbf{V}_3 \equiv \eE[\textbf{v}_3 \textbf{v}_3^\top]=\textbf{S}_{1,{\cal A}}+\textbf{S}_{2,{\cal A}},
 $$ 
with matrix $\textbf{S}_1$,   $\textbf{S}_2$ defined in Theorem \ref{Theorem2.1 Tang}. \\

By the following theorem we show that $\widehat\eb_n$ enjoys the oracle properties, i.e. its sparsity and the asymptotic normality of $\widehat\eb_{n,{\cal A}}$.
\begin{theorem}
	\label{Theorem2.3 Tang}
Under assumptions (A1)-(A7), if $(\lambda_n)_{n \in \N}$ and $\gamma >0$ are such that  $n^{-1/2} \lambda_n =O_\PP(1)$ and $n^{(\gamma -1)/2} \lambda_n  \rightarrow \infty $, thus:\\
	  (i) $\widehat\eb_n$ is sparse: $\lim_{n \rightarrow \infty} \PP\big[ \widehat {\cal A}_n ={\cal A}\big]=1$. \\
	  (ii) $\widehat\eb_{n,{\cal A}}$ is asymptotically normal:
	  \[
	 \sqrt{n} \big( \widehat{\eb}_{n,{\cal A}} - \eb^0_{\cal A}\big) \overset{\cal L} {\underset{n \rightarrow \infty}{\longrightarrow}}  {\cal N}_{|{\cal A}|} \bigg(- \eE^{-1}_\varepsilon[h_\tau(\varepsilon)] l_0  {\eo^0_{\cal A}}^\top \sign(\eb^0_{\cal A})\eE_\eX[\eX_{{\cal A}}\eX^\top_{{\cal A}}]^{-1}, \eXL \bigg), 
	  \]
 with $l_0 \equiv \lim_{n \rightarrow \infty} n^{-1/2} \lambda_n$, the  $|{\cal A}|$-square matrix  $\eXL \equiv \eE^{-2}_\varepsilon[h_\tau(\varepsilon)]\eE_\eX[\eX_{{\cal A}}\eX^\top_{{\cal A}}]^{-1}  \textbf{V}_3 \eE_\eX[\eX_{{\cal A}}\eX^\top_{{\cal A}}]^{-1}$ and $\eo^0_{\cal A} \equiv \lim_{n \rightarrow \infty} \widehat \eo_{n,{\cal A}} =|\eb^0_{\cal A}|^{-\gamma}$.  
\end{theorem}
 Recall that the same assumptions on $\gamma$ and $\lambda_n$ were considered for a classical linear model estimated by the adaptive  LASSO expectile method in the works of \cite{Liao-Park-Choi.19} and \cite{Ciuperca.21}.
 \begin{remark1}
For the particular case $\tau=1/2$, $\gamma =1$, $l_0=0$, by Theorems \ref{Theorem2.2 Tang} and \ref{Theorem2.3 Tang} we obtain the results proved by \cite{Johnson.09} for the censored adaptive LASSO LS estimator. Moreover, always for $\gamma=1$ and $l_0=0$, similar results were obtained by \cite{Shows.Lu.Zhang.10} when $l_0 \geq 0$ and by \cite{Tang-Zhou-Wu.12} when $l_0 =0$ for the censored adaptive LASSO median and quantile estimators, respectively.   	
 \end{remark1}
From Theorem \ref{Theorem2.3 Tang}\textit{(ii)} we deduce that if $l_0 \neq 0$ then the estimator $\widehat{\eb}_{n,{\cal A}} $ of $\eb^0_{\cal A}$ is asymptotically biased. On the other hand, the variance-covariance matrix  $\eXL$ is the corresponding matrix obtained in Theorem \ref{Theorem2.1 Tang} for the estimator $\widetilde{\eb}_{n,{\cal A}} $, that is, $\eXL=\textbf{S}^{-1}_{3,{\cal A}}(\textbf{S}_{1,{\cal A}}+\textbf{S}_{2,{\cal A}})\textbf{S}^{-1}_{3,{\cal A}})$. Then, the  estimators  $\widetilde{\eb}_{n,{\cal A}} $ and  $\widehat{\eb}_{n,{\cal A}} $ have the same asymptotic variance matrix and   for $l_0=0$ they have the same asymptotic centered normal distribution. \\
Moreover, let us emphasize that an additional difficulty arises in the theoretical study of the estimators $\widetilde{\eb}_n$, $\widehat{\eb}_n$ due to the presence of the Kaplan-Meier estimator which depends on the censored observations.\\
The form of the asymptotic variance-covariance matrix  $\eXL$   of the estimator  $ \widehat{\eb}_{n,{\cal A}} $ and therefore that of $\widetilde{\eb}_n$ is complex, which can lead to difficulties in practical applications for hypothesis testing or constructing the confidence interval for $\eb$. There are several methods to estimate this matrix. Either use the bootstrap method as proposed in \cite{Shows.Lu.Zhang.10}, \cite{Chen.Jewell.05}, or use the jackknife method  as suggested in \cite{Wang.Ng.08}. Another possibility is to estimate it  by the plug-in method using the empirical sample averages for  $\eE_\varepsilon[h_\tau(\varepsilon)]$, $\eE_\eX[\eX_{{\cal A}}\eX^\top_{{\cal A}}]$, $\textbf{S}_{1,{\cal A}}$,  $\textbf{S}_{2,{\cal A}}$.\\

In addition to these theoretical results, the present paper is also motivated by a simulation study that confirms these results and which also shows its superiority over other censored estimation methods in the literature. This numerical study is presented in the following section. The practical interest of the censored adaptive LASSO expectile estimation method is supported by applications on real data in Section \ref{section_application}.
\section{Simulation studies}
\label{section_simus}
In this section, through Monte Carlo simulations,  we illustrate the theoretical properties obtained  for the estimators $\widetilde \eb_n$ and $\widehat \eb_n$ in Section \ref{section_results}. We also compare the performance of the censored adaptive LASSO expectile estimator with that of the censored adaptive LASSO quantile and censored adaptive LASSO LS estimators. The R software was used to conduct the simulations. \\
 In subsections \ref{subsect_ln}, \ref{subsect_error}, \ref{subsect_compar}, the numerical study concerns the censored adaptive LASSO expectile estimator $\widehat \eb_n$, while subsection \ref{subsect_estexp} concerns the censored expectile estimator $\widetilde \eb_n$. In subsection \ref{subsect_conclusions} we formulate the conclusions obtained after all these simulations.\\
 
In this section we consider the censored model:
\begin{equation}
	\label{eq1T}
	T^*_i = \beta^0_0+ \sum^p_{j=1} \beta^0_j X_{ji} +\varepsilon_i, \qquad i=1, \cdots, n,
\end{equation}
which will be estimated in two cases: non-zero intercept   ($\beta^0_0 \neq 0$) or zero intercept ($\beta^0_0=0$).\\
In the case $\beta^0_0=0$, the parameters $\beta^0_0, \beta^0_1, \cdots , \beta^0_p$ of (\ref{eq1T}) will be estimated assuming that we know a priori that the model does not have an intercept (called supposition without intercept) and also when we don't know, in which case, we leave the possibility in the   coefficient estimation  of the intercept estimation (called supposition with intercept).\\

For $i=1, \cdots , n$, we conduct simulations for the design $X_{ji} \sim {\cal N}(1,1)$ for any $j=1, \cdots p$ and the censoring variable ${\cal C}_i \sim {\cal U}[0,c_1]$, with the constant $c_1$ chosen to obtain an a priori fixed censoring rate. Unless otherwise specified, the true coefficients are:  $\beta^0_1=0.9$, $\beta^0_2=-2$, $\beta^0_3=0.5$, $\beta^0_4=1$, $\beta^0_5=-1$ and $\beta^0_j=0$ for any $j \in \{ 6, \cdots,  p\}$. For the model errors $(\varepsilon_i)_{1 \leqslant i \leqslant n}$ we consider the following two distributions: Uniform ${\cal U}[-1,2]- 1/6$  and standard Gumbel  ${\cal G}(0,1)$. These  two distributions   are  not centered and ${\cal G}(0,1)$ is asymmetric. Recall that $\eE[{\cal G}(0,1)]$ is equal to Euler's constant and  that if $\varepsilon_i \sim {\cal G}(0,1)$, then, conditioned by ${\cal C}_i$ and $\eX_i$,  the distribution of $Y_i$ is Weibull ${\cal W}(0,1)$. The values considered for the constant $c_1$ are such that the censoring rate is $10\%$ or $25\%$. For all these configurations we calculate $Y_i$ and $\delta_i$ for $i=1, \cdots , n$ to then estimate the parameter $\eb$. We compare the penalized   censored expectile method proposed in this paper (for the expectile index $\tau$,  which verifies assumption  (A7)) with the censored adaptive LASSO quantile method considered in  \cite{Tang-Zhou-Wu.12} but for a single value of the quantile index $\widetilde \tau$ such that $\eE[\e1_{\varepsilon < 0}]=\widetilde \tau$ and for the tuning parameter $n^{1/2-0.1}$ (unless otherwise stated). We also compare the censored adaptive LASSO expectile estimator with  the censored  adaptive  LASSO LS estimator which is in fact a special case of our method for $\tau =1/2$. Recall that \cite{Johnson.09} considered the censored adaptive LASSO LS estimator for the particular case $\gamma=1$ and $l_0=0$.\\
For each scenario, 100 Monte Carlo replications have been carried. \\

In all Figures from \ref{fig_ChA_ChAn_b0}(a) to \ref{fig_ChA_ChAn_L2b0}(d),  we have represented the results obtained for the censored adaptive LASSO expectile estimator with the symbol $\circ$, with the symbol \textcolor{red}{$\square$} for the censored adaptive LASSO LS estimator, and with \textcolor{green}{$\blacktriangle$} for the censored adaptive LASSO quantile estimator.   The tuning parameter  $\lambda_n $ and the power  $\gamma=2$ in the adaptive weights were chosen according to with the assumptions imposed in Theorem \ref{Theorem2.3 Tang}.\\
For Figures  \ref{fig_ChAAn12_Gumbel}(a) to \ref{fig_ChA_ChAn_L2b0}(d) we calculated on $M$ Monte Carlo replications:
\begin{itemize}
\item \textit{the percentage of true zeros}, computed by: 
\[
100 \frac{1}{M}\sum^M_{l=1} \frac{\big|\widehat {\cal A}^c_{n,l} \cap {\cal A}^c\big|}{| {\cal A}^c|},
\]
\item \textit{the percentage of false zeros}, computed by:
\[
100 \frac{1}{M}\sum^M_{l=1} \frac{\big|\widehat {\cal A}^c_{n,l} \cap {\cal A}\big|}{| {\cal A}|},
\]
where $\widehat{\cal A}^c_{n,l}$ represents the complementary of  set $\widehat{\cal A}_{n,l}$ obtained for the $l$th Monte Carlo replication. Then, $\widehat{\cal A}^c_{n,l}$ contains the indexes of the zero components of the  estimation $\widehat \eb_n$ obtained to the $l$th Monte Carlo replication.
\end{itemize}
 In Figures \ref{fig_ChAAn12_Gumbel} to \ref{fig_ChA_ChAn_L2b0}, the horizontal line $95$ has been drawn for the rate of true zeros and the line $5$ for the figures of false zeros. 
 A perfect estimation method would produce a value of 100 for the  true zeros percentage and a value of 0 for the false zeros percentage.\\
 
 In Figures \ref{fig_ChAAn12_Gumbel}, \ref{fig_ChAn_NULG}, \ref{fig_ChA_GGGG_ai}, \ref{fig_ChAAn_GGr_ai},  \ref{fig_ChA_evolcen} and \ref{fig_mean_sd_Gumbe1l0_r10} we consider that true model (\ref{eq1T}) is without  intercept, i.e. $\beta^0_0=0$, but the model will be estimated in two cases: assuming that there is no intercept and assuming that  it is possible that there is an intercept. Moreover, for  Figures \ref{fig_ChA_ChAn_b0} and \ref{fig_ChA_ChAn_L2b0} we take for model (\ref{eq1T}) that the intercept is $\beta^0_0=2$.  For Figures  \ref{fig_ChAAn12_Gumbel} to  \ref{fig_ChA_evolcen} and Table \ref{Tab2}, we consider the index set ${\cal A}=\{1, 2, 3, 4, 5\}$. 
 \subsection{Numerical study on the choice of the tuning parameter $\lambda_n$}
 \label{subsect_ln}
\begin{figure}[h!] 
	\begin{tabular}{cc}
		\includegraphics[width=0.45\linewidth,height=4.5cm]{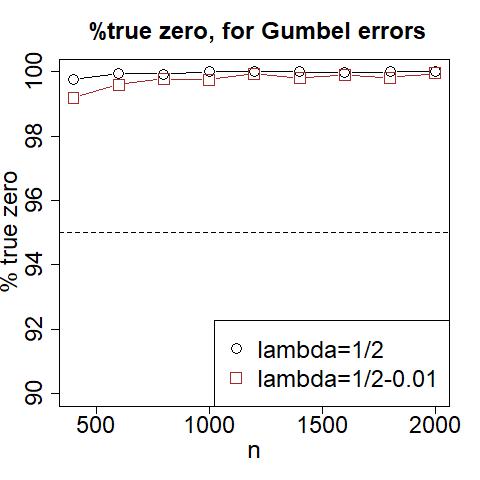} &
		\includegraphics[width=0.45\linewidth,height=4.5cm]{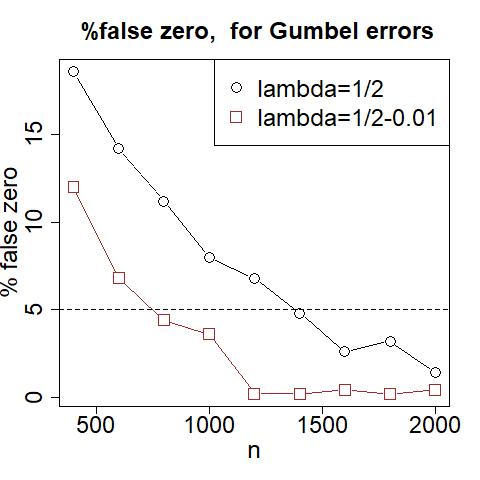} \\
		{\small (a) $\%$ of true zeros, supposition with intercept.} &
		{\small (b) $\%$ of false zeros, supposition with intercept.}\\
		& \\
		\includegraphics[width=0.45\linewidth,height=4.5cm]{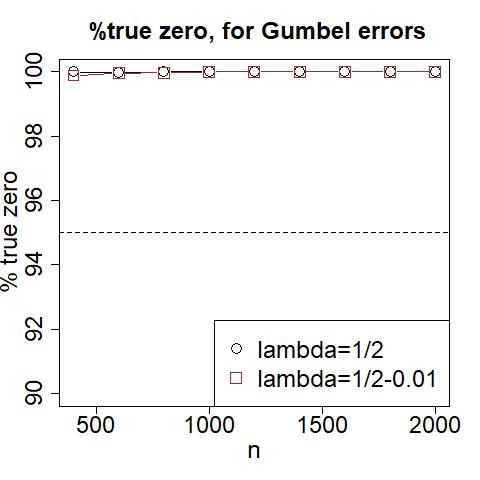} &
		\includegraphics[width=0.45\linewidth,height=4.5cm]{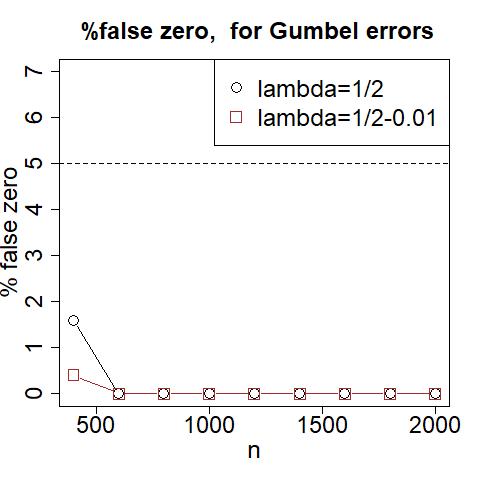} \\
		{\small (c) $\%$ of true zeros, supposition without intercept.} &
		{\small (d) $\%$ of false zeros, supposition without intercept.}
	\end{tabular}
	\caption{\small Percentage evolution of the true and false zeros with respect to $n$ for two sequences  $\lambda_n$  ($\circ$ for $\lambda_n=n^{1/2}$, \textcolor{red}{$\square$} for $\lambda_n =n^{1/2- 0.01}$)     by censored adaptive LASSO expectile method, for model without intercept ($\beta^0_0=0$),  when $p=50$, censoring rate $25\%$, $\varepsilon \sim{\cal G}(0,1)$.	}
	\label{fig_ChAAn12_Gumbel}
\end{figure}

\begin{figure}[h!] 
 	\begin{tabular}{cc}
 		\includegraphics[width=0.45\linewidth,height=4.5cm]{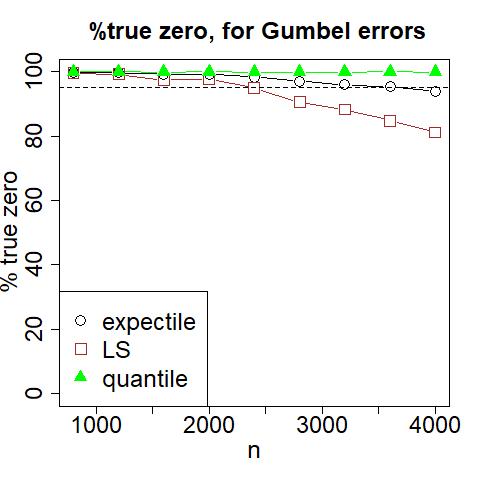} &
 		\includegraphics[width=0.45\linewidth,height=4.5cm]{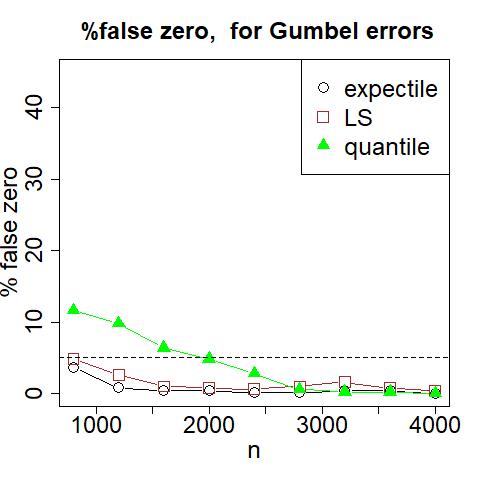} \\
 		{\small (a) $\%$ of true zeros for $\lambda_n=n^{1/2-0.1}$.} &
 		{\small (b) $\%$ of false zeros for $\lambda_n=n^{1/2-0.1}$.}\\
 		& \\
 		\includegraphics[width=0.45\linewidth,height=4.5cm]{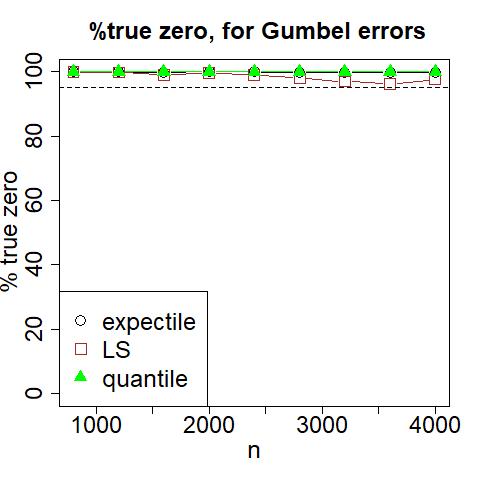} &
 		\includegraphics[width=0.45\linewidth,height=4.5cm]{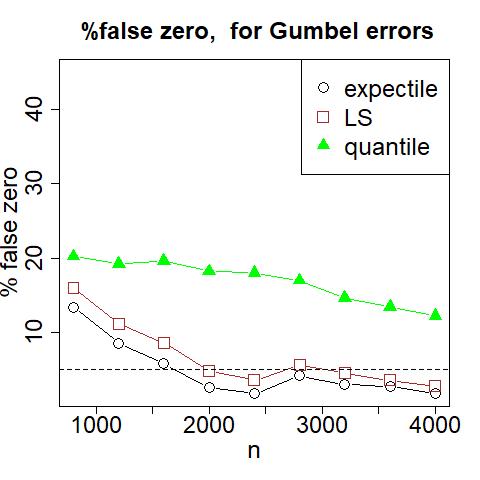} \\
 		{\small (c) $\%$ of true zeros for $\lambda_n=n^{1/2}$.} &
 		{\small (d) $\%$ of false zeros for $\lambda_n=n^{1/2}$.}
 	\end{tabular}
 	\caption{\small Percentage evolution of the true and false zeros with respect to $n$ for two sequences  $\lambda_n$ by three censored adaptive LASSO estimation methods,  when $\varepsilon \sim{\cal G}(0,1)$, $p=50$,  model with intercept ($\beta_0^0=2$) and censoring rate  $25\%$.	}
 	\label{fig_ChA_ChAn_b0}
 \end{figure}

 Based on the assumptions and  the asymptotic result of Theorem \ref{Theorem2.3 Tang}, in order to study the choice of the tuning parameter $\lambda_n$  on the automatic variable selection, we consider $\lambda_n =n^{1/2}$, i.e. $l_0=1$ and $\lambda_n=n^{1/2-0.1}$, i.e. $l_0=0$. From Figures \ref{fig_ChAAn12_Gumbel}(c) and \ref{fig_ChAAn12_Gumbel}(d) we deduce that the estimation of zero and non-zero coefficients is practically identical for the two  tuning parameter sequences when the assumed model is without intercept.  Thus,  if $\beta^0_0 =0$ and we force that the estimated model is without intercept, then  we can choose any  $\lambda_n$ with condition that $n^{-1/2}\lambda_n=O(1)$. On the other hand, if during the estimation we leave the possibility of an intercept, then for $l_0=0$, the rate $5\%$ of false zeros is achieved for a $n$ smaller than for $l_0=1$. In other words, convergence towards a false zero rate of $5\%$ is slower when $l_0 >0$. The detection of the true zeros is not disturbed by the choice of $\lambda_n$ (Figures \ref{fig_ChAAn12_Gumbel}(a) and \ref{fig_ChAAn12_Gumbel}(b)). Hence, the results of Figures \ref{fig_ChAAn12_Gumbel}(a)-\ref{fig_ChAAn12_Gumbel}(d)  confirm the statements announced by Theorem \ref{Theorem2.3 Tang}(i). These results are complemented by those of Table \ref{Tab2} where, in order to evaluate the choice influence of  $\lambda_n$ on the precision of the censored adaptive LASSO expectile estimation, we calculate $\| \widehat \eb_n -\eb^0\|$ and the standard-deviation of  $(\widehat{\eb}_n -\ebo)_{\cal A}$. First of all, from this table we obtain  confirmation of  convergence of  $\widehat\eb_n$ towards $\ebo$ stated by Theorem \ref{Theorem2.2 Tang} and the asymptotic bias of the estimator $\widehat \eb_{n,{\cal A}}$ stated by Theorem \ref{Theorem2.3 Tang}(ii). On the other hand, the standard deviation is the same for the two sequences considered for $\lambda_n$, which is consistent with the result of Theorem \ref{Theorem2.3 Tang}(ii).\\
 Figures  \ref{fig_ChA_ChAn_b0}(a) and \ref{fig_ChA_ChAn_b0}(d)  display the percentage of  true and false zeros by the three censored adaptive LASSO estimation  methods based on two tuning sequences $\lambda_n$ (inclusively for the quantile method)    $\lambda_n = n^{1/2}$ and   $\lambda_n = n^{1/2-0.1}$  when model (\ref{eq1T}) has an intercept, $\beta^0_0=2$. 
For $\lambda_n = n^{1/2}$, then all three methods detect over $95\%$ of the  true zeros, while for  $\lambda_n = n^{1/2-0.1}$, as $n \geq 2500$,  the true zero detection rate of the censored adaptive LASSO LS method decreases as $n$ increases and this rate is less than $95\%$ (see Figures \ref{fig_ChA_ChAn_b0}(a) and \ref{fig_ChA_ChAn_b0}(c)). Regarding the percentages of false zeros (Figures \ref{fig_ChA_ChAn_b0}(b) and \ref{fig_ChA_ChAn_b0}(d)), by the censored adaptive LASSO quantile method, we obtain  when $\lambda_n=n^{1/2}$ that this rate is greater than $5\%$ and less than $5\%$ when $\lambda_n=n^{1/2-0.1}$, if $n >2000$. This confirms the assumptions on $\lambda_n$ considered by \cite{Tang-Zhou-Wu.12}. By the censored adaptive LASSO expectile and LS methods, the percentage of false zeros is less than $5\%$  when $\lambda_n=n^{1/2-0.1}$  for any $n$ (Figure \ref{fig_ChA_ChAn_b0}(b))  and when $\lambda_n=n^{1/2}$ for $n \geq 1600$ by expectile, for $n \geq 2000$ by LS   (Figure \ref{fig_ChA_ChAn_b0}(d)). Therefore, in order to correctly choose, for any value of $n$, the true zeros and non-zero coefficients of   (\ref{eq1T}) when the model contains intercept, it is  best  to take  $\lambda_n = n^{1/2-0.1}$. In this case the expectile method gives excellent results. When $\lambda_n = n^{1/2-0.1}$, the quantile technique detects more than $5\%$ of non-zero coefficients   as zero for $n <2000$, while by the LS technique, we obtain that more than   $5\%$  zeros are estimated as non-zero for $n > 2500$. When $\lambda_n = n^{1/2}$ and $n \geq 1500$, the penalized expectile technique gives very good results which are better than those obtained by the LS and quantile techniques.\\
Following the conclusions of this numerical study, in all the simulations that follow and in the applications of Section \ref{section_application}  we will take the tuning parameter $\lambda_n=n^{1/2-0.1}$ for the penalties of the three loss functions. 
 \begin{table}
 	\begin{center}
 		\begin{tabular}{|c|cc|cc|}\hline  
 			$n$	&  \multicolumn{2}{c|}{$\lambda_n=n^{1/2} $}  & \multicolumn{2}{c|}{$\lambda_n=n^{1/2-0.1}$}  \\ 
 			\cline{2-5}
 			&   \textit{L2}   &\textit{SD} & \textit{L2} & \textit{SD}    \\ \hline 
 			400 & 0.34 & 0.08 & 0.28 & 0.07   \\
 			1000 & 0.24 & 0.05 & 0.20 & 0.05  \\
 			2000 & 0.25 & 0.10 & 0.23 & 0.10 \\  \hline
 		\end{tabular} 
 	\end{center}
 	\caption{\small Accuracy results  ($L2\equiv \| \widehat \eb_n -\eb^0\|$, $SD \equiv sd((\widehat{\eb}_n -\ebo)_{\cal A})$)  of   $  \widehat \eb_n  $  obtained by 100 Monte Carlo replications  for censored adaptive LASSO expectile method, when $\varepsilon \sim {\cal G}(0,1)$, $p=50$, $\beta^0_0=0$, ${\cal A}=\{1, 2, 3, 4, 5\}$, censoring rate $25\%$,  models estimated without intercept.}
 	\label{Tab2} 
 \end{table}

 \subsection{Numerical study with respect to model error distributions}
 \label{subsect_error}
\begin{figure}[h!] 
	\begin{tabular}{cc}
		\includegraphics[width=0.45\linewidth,height=4.5cm]{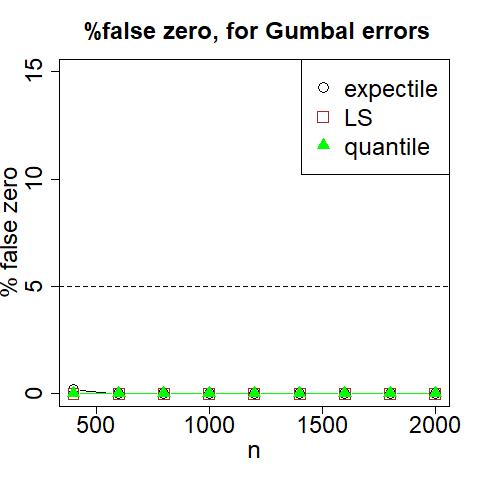} &
		\includegraphics[width=0.45\linewidth,height=4.5cm]{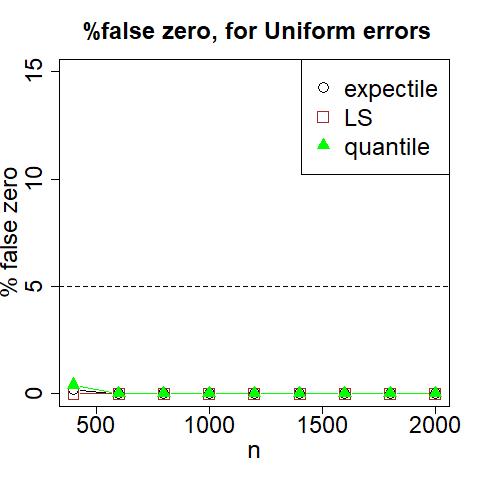} \\
		{\small (a) When	$\varepsilon \sim{\cal G}(0,1) $, supposition without intercept.} &
		{\small (b) When $\varepsilon \sim {\cal U}[-1,2]-1/6$, supposition without intercept.} \\
		& \\
			\includegraphics[width=0.45\linewidth,height=4.5cm]{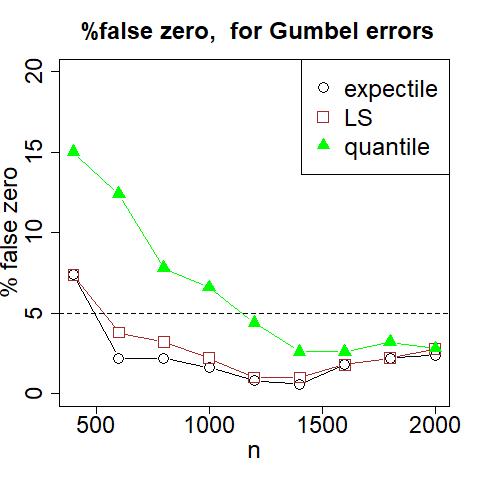} &
		\includegraphics[width=0.45\linewidth,height=4.5cm]{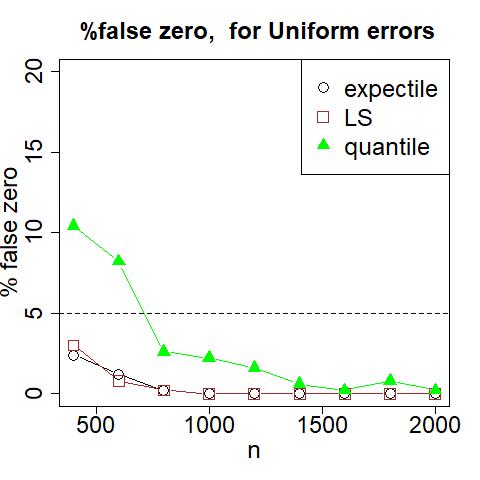} \\
		{\small (c) When	$\varepsilon \sim{\cal G}(0,1) $, supposition with intercept.} &
		{\small (d) When  $\varepsilon \sim {\cal U}[-1,2]-1/6$, supposition with intercept.}
	\end{tabular}
	\caption{\small Percentage evolution of the false zeros  by three censored adaptive LASSO estimation methods, for model without intercept ($\beta^0_0=0$), when $p=50$ and censoring rate is $25\%$.	}
	\label{fig_ChAn_NULG}
\end{figure}

In Figures  \ref{fig_ChAn_NULG}(a) - \ref{fig_ChAn_NULG}(d) we represent with respect to $n \in \{ 400, 600, \cdots, 2000\}$ the percentage of false zeros  when $p=50$, censoring  rate   equal to $25\%$, model errors ${\cal G}(0,1)$ and ${\cal U}[-1,2]-1/6$, without intercept $\beta_0^0=0$, but  estimated model supposed without and with intercept.
The corresponding figures for  the percentage of  true zeros for the same configurations  are not shown  because this percentage  is always $100\%$ for all values of $n$.  Regarding the detection of false zeros, if we assume that the estimated model is without intercept, this rate is 0 or very close to it (Figures  \ref{fig_ChAn_NULG}(a) and  \ref{fig_ChAn_NULG}(b)). On the other hand, by offering it the possibility of having an intercept (Figures  \ref{fig_ChAn_NULG}(c) and  \ref{fig_ChAn_NULG}(d)), the censored adaptive LASSO LS and expectile methods give similar results, which are better than those by the censored adaptive LASSO quantile method which makes more false zero detections, especially when $n<1500$.
\subsection{Comparative numerical study by varying $p$, censoring rate, $\beta^0_0$, $n$}
\label{subsect_compar}
\begin{figure}[h!] 
 	\begin{tabular}{cc}
 			\includegraphics[width=0.45\linewidth,height=4.5cm]{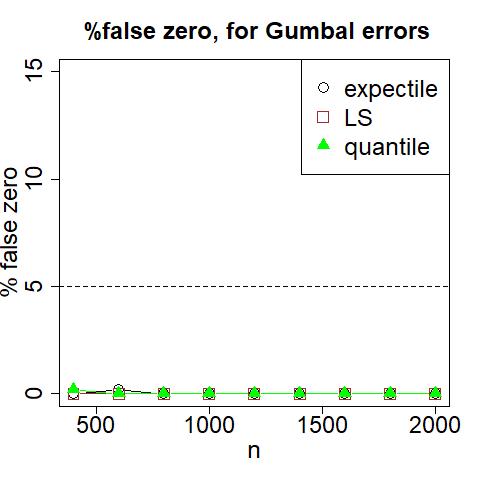} &
 		\includegraphics[width=0.45\linewidth,height=4.5cm]{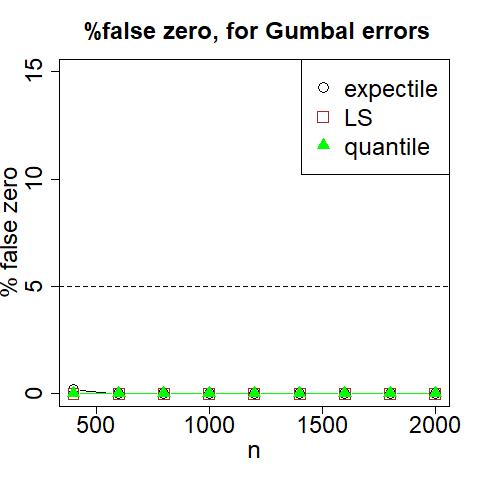} \\
 		{\small (a)  Supposition without intercept, $p=50$.} &
 		{\small (b) Supposition without intercept, $p=150$.} \\
 		& \\
 		\includegraphics[width=0.45\linewidth,height=4.5cm]{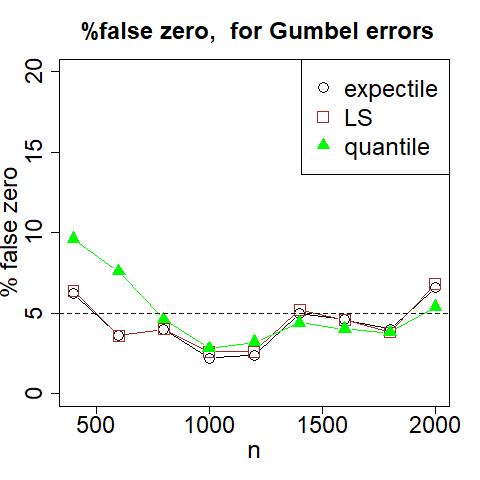} &
 		\includegraphics[width=0.45\linewidth,height=4.5cm]{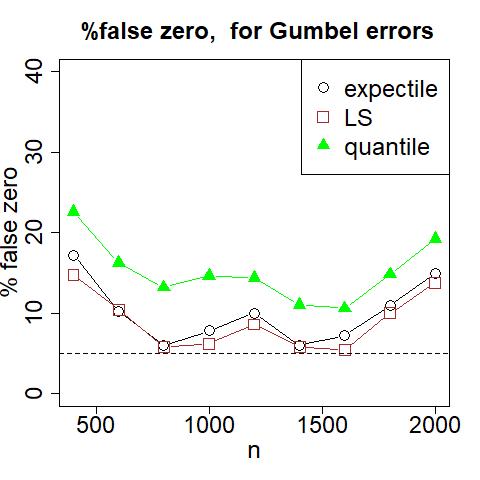} \\
 		{\small (c) Supposition with intercept, $p=50$.} &
 		{\small (d) Supposition with intercept, $p=150$.}
 	\end{tabular}
 	\caption{\small Percentage evolution of the false zeros  by three censored adaptive LASSO estimation methods, for model without intercept ($\beta^0_0=0$),  when $\varepsilon \sim{\cal G}(0,1)$ and censoring rate is  $10\%$.	}
 	\label{fig_ChA_GGGG_ai}
 \end{figure}

\begin{figure}[h!] 
 	\begin{tabular}{cc}
 		\includegraphics[width=0.45\linewidth,height=4.5cm]{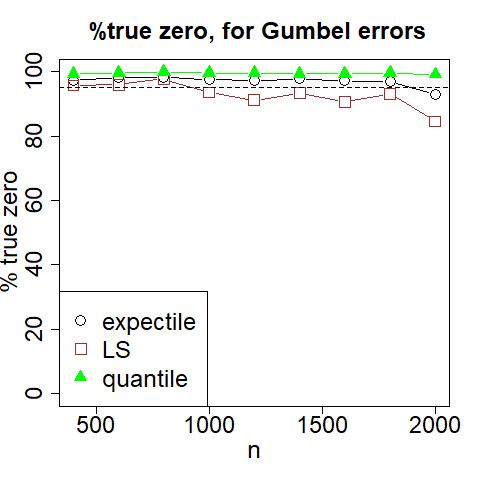} &
 		\includegraphics[width=0.45\linewidth,height=4.5cm]{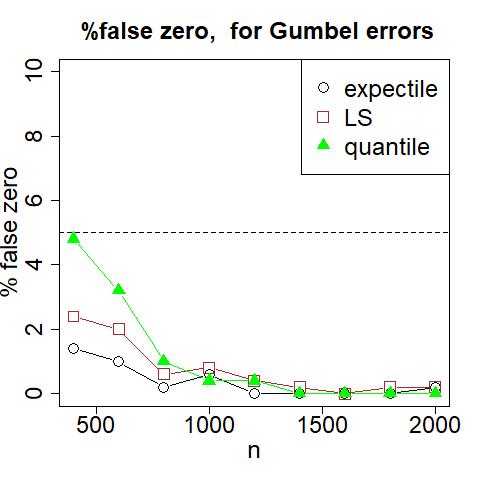} \\
 		{\small (a)  $\%$ of true zeros for   censoring rate $25\%$.} &
 		{\small (b) $\%$ of false zeros for  censoring rate $25\%$.} \\
 		& \\
 		\includegraphics[width=0.45\linewidth,height=4.5cm]{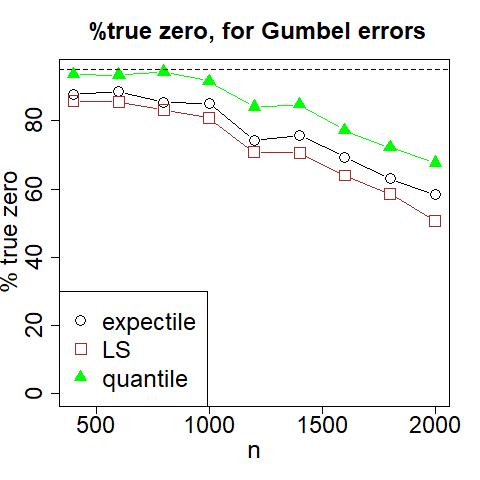} &
 		\includegraphics[width=0.45\linewidth,height=4.5cm]{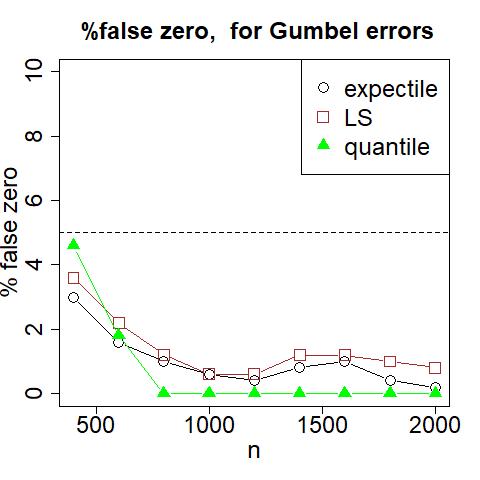} \\
 		{\small (c) $\%$ of true zeros for    censoring rate $10\%$.} &
 		{\small (d) $\%$ of false zeros for   censoring rate $10\%$.}
 	\end{tabular}
 	\caption{\small Percentage evolution  of the true and false zeros  by three censored adaptive LASSO estimation methods, for model without intercept ($\beta^0_0=0$), supposition with intercept, when $\varepsilon \sim{\cal G}(0,1) $ and $p=10$.	}
 	\label{fig_ChAAn_GGr_ai}
 \end{figure}

Starting with Figure   \ref{fig_ChA_GGGG_ai} we will focus on the case $\varepsilon \sim  {\cal G}(0,1)$  because this is the most common case  in censorship models. We will vary the number of explanatory variables of the model, but the  non-zero coefficients will always be the first five. We also vary the censoring rate, considering as values:  $10\%$ or $25\%$. \\
In Figures  \ref{fig_ChA_GGGG_ai}(a) -  \ref{fig_ChA_GGGG_ai}(b) we represent the percentage  of  false zeros when $p \in \{50, 150\}$, $\varepsilon \sim {\cal G}(0,1) $,  censoring rate equal to $10\%$, $\beta^0_0=0$, and when we suppose that the estimated model is without intercept. In Figures  \ref{fig_ChA_GGGG_ai}(c) - \ref{fig_ChA_GGGG_ai}(d) we have  the percentage of false zeros  for the same configurations but leaving the possibility  that the estimated model has an intercept.  To investigate the effect of the number of zero components of  $\ebo$ and the effect of the censoring rate on the sparsity of the three censored adaptive LASSO estimators, these results should be compared with those of Figures  \ref{fig_ChAAn_GGr_ai}(a)-\ref{fig_ChAAn_GGr_ai}(d) where $p=10$.\\
If the censorship rate is $10\%$ and $p$ is either 50 or 150, all three estimation methods detect over $95\%$ of true zeros (figures not shown) and have less than $1\%$ false zeros when assuming no intercept in the estimation model (Figures \ref{fig_ChA_GGGG_ai}(a) -  \ref{fig_ChA_GGGG_ai}(b) ). If we leave the possibility of intercept (which is not present in the true model), we can see from Figures \ref{fig_ChA_GGGG_ai}(c) and \ref{fig_ChA_GGGG_ai}(d) we deduce that the censored adaptive LASSO LS and expectile methods give similar results, results which are better than by censored adaptive LASSO quantile method. On the other hand,  the methods detect over $95\%$ of the real zeros and therefore we do not display the figures. \\
When $p=10$, i.e. there are few real zeros, we study the percentage evolution  of  true and false zeros  for three censored adaptive LASSO estimation methods.  This is performed assuming the model lacks an intercept, when  $\varepsilon  \sim {\cal G}(0,1) $, and the censoring rate is either $25\%$ or $10\%$. Using the three methods, we  identify at least  $95\%$ of true zeros and commit less than $1\%$ of false zeros (figures not shown).  The results deteriorate slightly if we give the possibility of intercept (Figures \ref{fig_ChAAn_GGr_ai}(a)-\ref{fig_ChAAn_GGr_ai}(d)).\\

\begin{figure}[h!] 
	\begin{tabular}{cc}
		\includegraphics[width=0.45\linewidth,height=4.5cm]{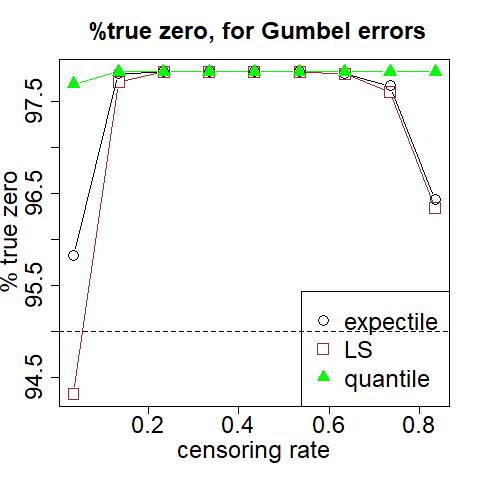} &
		\includegraphics[width=0.45\linewidth,height=4.5cm]{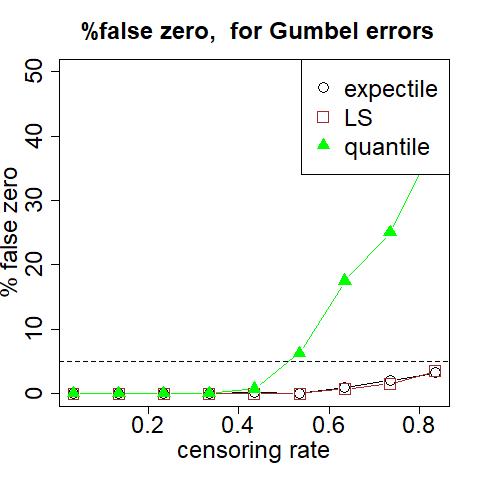} \\
		{\small (a) $\%$ of true zeros,  supposition without intercept.} &
		{\small (b) $\%$ of false zeros,  supposition without intercept.} \\
		& \\
		\includegraphics[width=0.45\linewidth,height=4.5cm]{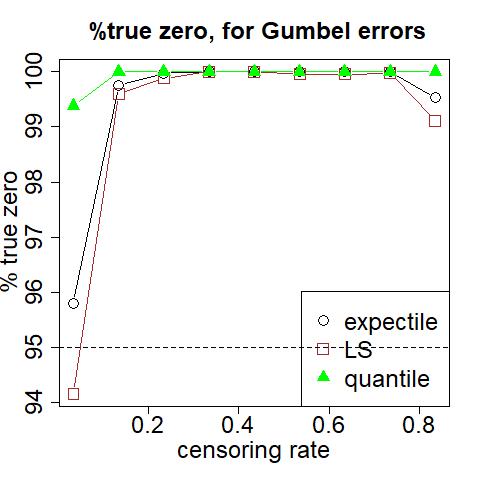} &
		\includegraphics[width=0.45\linewidth,height=4.5cm]{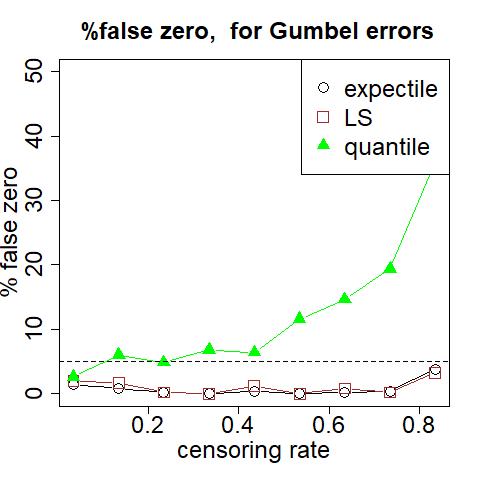} \\
		{\small (c) $\%$ of true zeros, supposition with intercept.} &
		{\small (d) $\%$ of false zeros, supposition with intercept.}
	\end{tabular}
	\caption{\small Percentage evolution  with respect to the censoring rate of the true and false zeros  by three censored adaptive LASSO estimation methods, when $\varepsilon \sim{\cal G}(0,1)$, $n=1000$, $p=50$, model without intercept ($\beta^0_0=0$).	}
	\label{fig_ChA_evolcen}
\end{figure}

\begin{figure}[h!] 
	\begin{tabular}{cc}
		\includegraphics[width=0.45\linewidth,height=4.5cm]{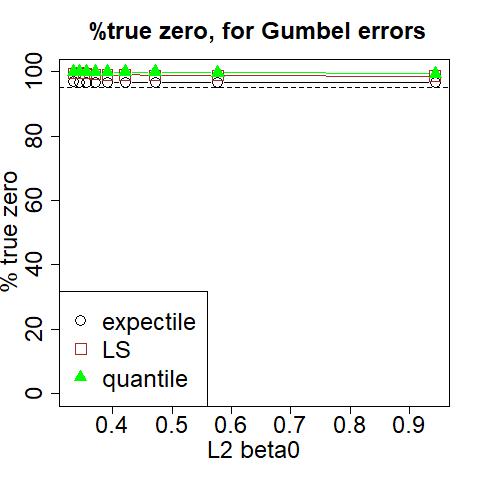} &
		\includegraphics[width=0.45\linewidth,height=4.5cm]{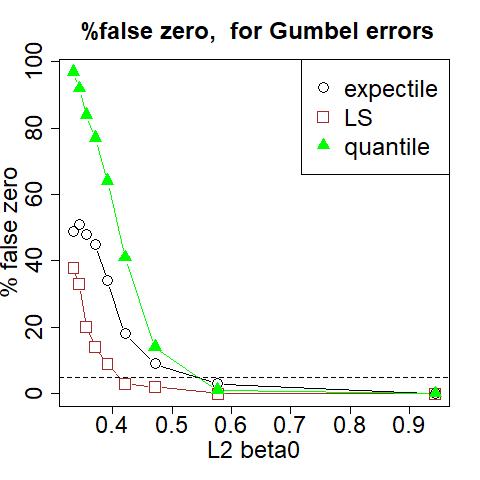} \\
		{\small (a) $\%$ of true zeros.} &
		{\small (b) $\%$ of false zeros.}
		\end{tabular}
	\caption{\small Percentage evolution    of the true and  false zeros with respect to L2 beta0 $\equiv \| \ebo\|$  by three censored adaptive LASSO estimation methods,  when $\varepsilon \sim{\cal G}(0,1)$, $p=50$, ${\cal A}=\{1,2\}$,  model with intercept ($\beta^0_0=2$) and censoring rate is  $25\%$.	}
	\label{fig_ChA_ChAn_L2b0}
\end{figure}

From Figures  \ref{fig_ChA_evolcen}(a) and \ref{fig_ChA_evolcen}(c)  we deduce that for all three methods, the detection of true zeros does not evolve with the censoring rate. Moreover, let's make a very important remark that the percentage of false zeros  for the censored adaptive LASSO expectile and LS estimators  is less than $5\%$ for any value of the censoring rate (Figures of \ref{fig_ChA_evolcen}(b) and \ref{fig_ChA_evolcen}(d)).  On the other hand, the censored adaptive LASSO quantile method detects more and more false zeros as the censoring rate increases.\\
Let's also study these methods with respect to the value of  $\|\eb^0\|$. For this we take in model (\ref{eq1T})  a single non-zero coefficient, more precisely ${\cal A}=\{1\}$. The considered model contains the intercept $\beta^0_0=2$, the value of the non-zero coefficient is $\beta_1^0=1/\log k - 1/k$, the tuning parameter is $\lambda_n=(150k)^{1/2-0.1}$, for $k \in \{2, 3 \cdots, 10\}$, $n=150k$ and $\beta_j^0=0$ for all $j \in \{2, \cdots, 50\}$. Then the values of $\beta_1^0$ are between $0.33$ and $0.94$. Once again, we observe that the true zero rate exceeds $95\%$ for the three estimation methods (Figure \ref{fig_ChA_ChAn_L2b0}(a)). For low values of $\|\ebo\|$, the coefficient of the variable $X_1$ can be  shrunk to 0 by the three estimation methods, with the rate of  false zeros in descending order: LS, expectile, and quantile methods. The percentage of false zeros decreases toward 0 as $\|\ebo\|$ increases  (Figure \ref{fig_ChA_ChAn_L2b0}(b)).  
 
 \begin{figure}[h!] 
	\begin{tabular}{cc}
		\includegraphics[width=0.45\linewidth,height=4.5cm]{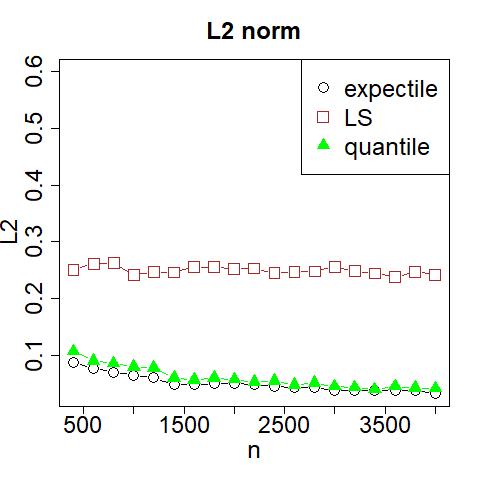} &
		\includegraphics[width=0.45\linewidth,height=4.5cm]{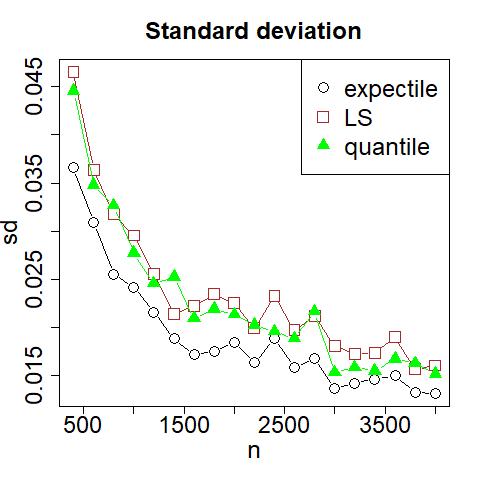} \\
		{\small (a) Evolution of  $L2\equiv \| \widetilde \eb_n -\ebo\|$.} &
		{\small (b) Evolution of $sd \equiv sd( \widetilde \eb_n -\ebo)$.}
	\end{tabular}
	\caption{\small Evolution of the accuracy of the parameter estimations  by three censored unpenalized estimation methods, when $\varepsilon \sim{\cal G}(0,1)$, $p=2$,  $\beta^0_0=0$,  censoring rate $10\%$.	}
	\label{fig_mean_sd_Gumbe1l0_r10}
\end{figure}
\subsection{Numerical study of $\widetilde \eb_n$}
\label{subsect_estexp}
We now study the evolution with $n$ of the  accuracy  of the censored expectile estimator $\widetilde \eb_n$ for   Gumbel error ${\cal G}(0,1) $, $p=2$: $\beta^0_1=5\log(n)$, $\beta^0_1=\log(n)$, $X_{1i} \sim {\cal N}(1,1)$,  $X_{2i} \sim {\cal N}(1,5)$, $i=1, \cdots , n$. The considered values for $\ebo$ allow to investigate the value effect of the parameter norm on the estimator accuracy. Furthermore, for 100 Monte Carlo replications, we compare the censored expectile estimator with the censored least squares and censored quantile estimators by calculating  the Euclidean norm $\|\widetilde{\eb}_n -\ebo\|$ and the standard deviation of $\widetilde{\eb}_n -\ebo$. From the results reported in  Table \ref{Tab1} and Figure \ref{fig_mean_sd_Gumbe1l0_r10} (symbol $\circ$ for  the censored  expectile estimator, \textcolor{red}{$\square$} for the  censored LS estimator, \textcolor{green}{$\blacktriangle$} for the  censored quantile estimator), we deduce that the  censored expectile estimators are more accurate than the other two estimators, especially for small values of $n$. The evolution of $\| \widetilde \eb_n -\ebo \| $ with $n$ in Table \ref{Tab1}  also supports the consistency of the estimator $\widetilde \eb_n$.  
\begin{table}
	\begin{center}
		\begin{tabular}{|c|ccc|ccc|}\hline  
			$n$	&  \multicolumn{3}{c|}{$\| \widetilde \eb_n -\ebo \| $}  & \multicolumn{3}{c|}{$sd( \widetilde \eb_n -\ebo )$}  \\ 
			\cline{2-7}
			& \textit{expectile} & \textit{LS} & \textit{quantile} &\textit{expectile} & \textit{LS} & \textit{quantile}    \\ \hline 
			10 & 0.59 & 0.70 & 0.68 & 0.21 & 0.23 &  0.23 \\
			50 & 0.25 & 0.37 & 0.30 & 0.09 & 0.11 &  0.12 \\
			100 & 0.16 & 0.29 & 0.19 & 0.06 & 0.08 &  0.08 \\
			200 & 0.12 & 0.27 & 0.14 & 0.05 &  0.06 &  0.05 \\ \hline
		\end{tabular} 
	\end{center}
	\caption{\small Accuracy results   of   $  \widetilde \eb_n  $  obtained by 100 Monte Carlo replications  for three censored unpenalized estimation methods, for model without intercept ($\beta^0_0=0$),  when $\varepsilon \sim {\cal G}(0,1)$, $p=2$, censoring rate $25\%$.}
	\label{Tab1} 
\end{table}

\subsection{Conclusion of simulations}
\label{subsect_conclusions}
The simulation results confirm the consistency and  sparsity of the censored adaptive LASSO expectile estimator $ \widehat \eb_n$. The consistency of the censored expectile estimator $\widetilde \eb_n$ is also shown.\\
The censored adaptive LASSO expectile estimator produces fewer false zeros and detects the zero coefficients better than the estimators corresponding to the quantile and LS methods. These detections do not depend on the distribution of the model error $\varepsilon$  for each of the three adaptive LASSO methods,  but the quantile method produces more false zeros than the  expectile and LS methods. On the other hand, for a given censoring rate,  when $\varepsilon \sim{\cal G}(0,1) $ and $\lambda_n$, $|{\cal A}|$ are fixed, then the percentage of false zeros does not depend on $p$ and $n$. We also obtained that by the censored adaptive LASSO expectile method,  the percentage of true zeros and that of  false zeros   do not evolve with the censoring rate when  $|{\cal A}|$, $n$, $p$ are fixed and $\varepsilon \sim{\cal G}(0,1) $. Let us emphasize that the censored adaptive LASSO quantile method detects more and more false zeros as the censoring rate increases. \\
For the choice of the tuning parameter sequence $(\lambda_n)_{n \in \N}$ in relation (\ref{ebn}) we recommend $\lambda_n =o(n^{1/2})$. \\
We conclude by emphasizing that  the censored expectile estimator is more accurate than the censored quantile and censored LS estimators.
\section{Applications on real data}
\label{section_application}
In order to support the practical interest of our method,  in this section we present the application of the censored adaptive LASSO expectile method on three survival databases.  The results are compared to those obtained by the  censored adaptive LASSO quantile and LS methods.\\
Note that for the following three applications, only continuous explanatory variables  are considered. These variables are standardized. 
\subsection{Primary Biliary Cholangitis Data}
We will use the non-penalized censored  expectile method and the one censored adaptive LASSO expectile on the dataset \textit{pbc} from the R package  \textit{survival}. The censoring variable ${\cal C}$ is in this case the variable \textit{status} which indicates the status at the endpoint  and which takes three values: 0/1/2 for censored, transplanted, dead, respectively. Initially, there are 418 patients in the study. If the patient received a transplant, then he is also considered as censored. Since the measured variables contain missing values, then the complete database  contains $n=276$ observations of which 165 are censored ($\delta_i =0$).  The failure time $T_i$ is the variable \textit{time} which gives the number of days between registration and the earlier of death. The  nine continuous explanatory variables of  censored model  (\ref{eq1T}) are: \textit{age, albumin, alk, ast, bili, chol, copper, platelet, protime}. See  the R package \textit{survival}  and its references for a more detailed description of the data.\\
For these nine variables, the estimation $\widehat \eb_n$ has seven zero components, more precisely only the variables \textit{albumin} and \textit{protime} have non-zero coefficients.  Considering censored model (\ref{eq1T}) only for these two explanatory variables and estimated by the censored expectile method, we obtain the estimations $2.314$ and $2.084$ of the coefficients, respectively.\\
Note also that the censored adaptive  LASSO  LS and quantile methods also select the explanatory variables \textit{albumin} and \textit{protime}.
\subsection{Myeloma}
In this subsection we consider two datasets from R packages on myeloma.
\subsubsection{Myeloma data of the R package "emplik"}
A multiple melanoma study is considered for 65 patients of whom 48 have died and 17 are alive at the end of the study. The survival time $T_i$ is measured in months.  The seven explanatory continuous variables of  model (\ref{eq1T}) are  \textit{AGE} (age at diagnosis) and expressions of six genes: \textit{LOGBUN} (log BUN at diagnosis), \textit{HGB} (hemoglobin at diagnosis),   \textit{LOGWBC} (log WBC at diagnosis), \textit{LOGPBM} (log percentage of plasma cells in bone marrow), \textit{PROTEIN} (proteinuria at diagnosis) and \textit{SCALC} (serum calcium at diagnosis).\\
 The nonzero censored adaptive LASSO expectile estimations are for the coefficients of the variables \textit{HGB, AGE, LOGWBC,  PROTEIN} and \textit{SCALC}. The coefficient estimations of these variables by the censored expectile method are: -1.3, 2.5, 16.9, 0.7, 2.9, respectively. \\
 Note that by the censored adaptive  LASSO LS  method the same five variables are selected, results confirmed by hypothesis tests, the residuals of the model being of normal distribution. By the censored adaptive LASSO quantile method only three variables are selected: \textit{LOGWBC,  PROTEIN} and \textit{SCALC}.
\subsubsection{Myeloma data of the R package "survminer"}
In this example, we study the survival time (in months) of patients with a multiple melanoma function of gene expression. There are 256 observations in total,  but we  only analyze 186 because the others have missing data. In the 186 observations, 51 patients have died and 135 are alive. The six continuous explanatory variables are: \textit{CCND1, CRIM1, DEPDC1, IRF4, TP53, WHSC1}.\\
By the censored adaptive LASSO expectile method, only the coefficient of {\it IRF4} is non-zero and its estimation by the censored expectile method is 0.49. Note that the same variable is selected by the censored adaptive LASSO LS method. On the other hand, the censored adaptive LASSO quantile method shrinks all coefficients to 0, meaning that no gene expression affects survival time.
\section{Proofs}
\label{section_proofs}
In this section we present the proofs of the results stated in Section \ref{section_results}.\\

\noindent {\bf Proof of  Theorem \ref{Theorem2.1 Tang}}
\textit{(i)} Let us consider the random $p$-vector $\widetilde \eu_n \equiv n^{1/2}(\widetilde \eb_n - \ebo)$. Thus, in order to prove the convergence rate of $\widetilde \eb_n $  we consider the parameter $\eb$ under the  form $\eb=\ebo+n^{-1/2} \eu$, with $\eu \in \R^p$ such that $\| \eu \|=c < \infty$. Then, let's study the following random process:
\[ 
	{\cal R}_n(\widehat G_n, \eb)  \equiv \sum^n_{i=1} \frac{\delta_i}{\widehat G_n (Y_i)} \bigg( \rho_\tau \big(\log(Y_i) -\eX^\top_i \eb\big) -\rho_\tau(\varepsilon_i)\bigg) 
 = \sum^n_{i=1} \frac{\delta_i}{\widehat G_n (Y_i)} \bigg( \rho_\tau \big( \varepsilon_i - n^{-1/2}\eX^\top_i  \eu\big) -\rho_\tau(\varepsilon_i)\bigg).
\]
We consider for a survival function $G$ and parameter $\eb$ the following random process:
\[
Q_n(G,\eb) \equiv \sum^n_{i=1} \frac{\delta_i}{G(Y_i)} \rho_\tau (\log(Y_i) -   \eX_i^\top \eb).
\]
Hence we have ${\cal R}_n(\widehat G_n, \eb)=Q_n(\widehat G_n,\eb)- Q_n(\widehat G_n,\ebo)$, which can be written:
\begin{align}
	\label{QQQ}
	Q_n(\widehat G_n,\ebo+n^{-1/2} \eu)- Q_n(\widehat G_n,\ebo) & =\left\{Q_n(G_0,\ebo+n^{-1/2} \eu)- Q_n(G_0,\ebo) \right\} \nonumber\\
	& \quad +\left\{Q_n(\widehat G_n,\ebo+n^{-1/2} \eu)- Q_n(G_0,\ebo+n^{-1/2} \eu) \right\} \nonumber\\
	& \quad - \left\{Q_n(\widehat G_n,\ebo)- Q_n(G_0,\ebo) \right\} \nonumber\\
	& \equiv {\cal Q}_{1n}+{\cal Q}_{2n}-{\cal Q}_{n}.
\end{align}
By elementary calculus, for $e \in \R$ and $t \rightarrow 0$, we have
\begin{equation}
	\label{GuZou}
	\rho_\tau(e -t)=\rho_\tau(e)+g_\tau(e)t+h_\tau(e) t^2/2+o(t^2).
\end{equation}
For term ${\cal Q}_{1n}$, written as:
\[
{\cal Q}_{1n}=\sum^n_{i=1} \frac{\delta_i}{G_0(Y_i)} \big( \rho_\tau (\varepsilon_i - n^{-1/2} \eX_i^\top \eu) - \rho_\tau(\varepsilon_i)\big),
\]
using the Cauchy-Schwarz inequality together with the fact that $\eE_{\cal C}[\delta_i]=G_0(Y_i)$, relation (\ref{GuZou}) and that $\|\eX\|$ is bounded by (A5), we get:
\begin{equation}
	\label{LGV}
{\cal Q}_{1n}= \sum^n_{i=1} \frac{\delta_i}{G_0(Y_i)} \bigg(  g_\tau(\varepsilon_i) \frac{\eX_i^\top \eu}{ \sqrt{n}} +\frac{1}{2} h_\tau(\varepsilon_i) \bigg(  \frac{\eX_i^\top \eu}{ \sqrt{n}} \bigg)^2+o_\PP (n^{-1})\bigg).
\end{equation}
On the other hand, using $\eE_{\cal C}[\delta_i]=G_0(Y_i)$ and $\eE_\varepsilon[g_\tau(\varepsilon_i)]=0$ of assumption (A7), we have:
\[
 \eE_\eX \bigg[\eE_\varepsilon \bigg[ \eE_{\cal C} \bigg[ \frac{\delta_i}{G_0(Y_i)} g_\tau (\varepsilon_i)  \eX_i^\top \eu   | \varepsilon_i, \eX_i \bigg] | \eX_i\bigg]\bigg]=\eE_\eX\big[ \eE_\varepsilon \big[ g_\tau (\varepsilon_i)\eX_i^\top \eu  |  \eX_i  \big]\big]=0.
 \]
We proceed in the same way for
\begin{align*}
\eE_\eX \bigg[\eE_\varepsilon \bigg[ \eE_{\cal C} \bigg[ \bigg( \frac{\delta_i}{G_0(Y_i)} g_\tau (\varepsilon_i)   \frac{\eX_i^\top \eu}{\sqrt{n}}\bigg)^2  | \varepsilon_i, \eX_i \bigg] | \eX_i\bigg]\bigg] & =\eE_\eX\bigg[ \eE_\varepsilon \bigg[  \frac{g^2_\tau (\varepsilon_i)}{G_0(Y_i)}   \frac{\eu \eX_i\eX_i^\top \eu}{n} |  \eX_i  \bigg]\bigg]\\
& =\eE_\varepsilon\big[ g^2_\tau(\varepsilon_i)\big]  \frac{\eu^\top}{n}  \eE_\eX \bigg[  \frac{\eX \eX^\top}{G_0(Y)}\bigg] \eu,
\end{align*}
 where we used that $\varepsilon_i$ is independent of $\eX_i$ of assumption (A1). Then, using the Central Limit Theorem (CLT), we get that
\[
\sum^n_{i=1} \frac{\delta_i}{G_0(Y_i)} g_\tau(\varepsilon_i) \frac{\eX_i^\top \eu}{ \sqrt{n}} =\textbf{W}^\top_1 \eu \big(1+o_\PP(1)\big),
\]
with the random $p$-vector $\textbf{W}_1 \sim {\cal N}_p(\oo_p, \textbf{S}_1)$.  \\
On the other hand, by the law of large numbers, we have:
\[
\frac{1}{n} \sum^n_{i=1} \frac{\delta_i}{G_0(Y_i)} h_\tau(\varepsilon_i) \eu^\top \eX_i \eX^\top_i \eu \overset{\PP} {\underset{n \rightarrow \infty}{\longrightarrow}} \eE\big[h_\tau(\varepsilon) \eu^\top \eX \eX^\top \eu \big]=\eE_\varepsilon \big[h_\tau(\varepsilon)\big] \eu^\top \eE_\eX\big[\eX \eX^\top\big] \eu.
\]
Taking into account these last relations together with relation (\ref{LGV}), then the term ${\cal Q}_{1n}$ becomes:
\begin{equation}
	\label{Q1n}
	{\cal Q}_{1n}=\big( \textbf{W}^\top_1 \eu+ \eE_\varepsilon \big[h_\tau(\varepsilon)\big] \eu^\top \eE_\eX\big[\eX \eX^\top\big] \eu\big) \big(1+o_\PP(1)\big).
\end{equation}
For terms ${\cal Q}_{2n}$ and ${\cal Q}_{3n}$, taking into account assumptions (A4) and (A6), we use the Taylor expansion of $\widehat G_n(Y_i)$ with respect to $G_0(Y_i)$:
\begin{equation}
\label{rell}
\begin{split}
	\sqrt{n} \bigg( \frac{1}{\widehat G_n(Y_i)} - \frac{1}{G_0(Y_i)}\bigg)& = - \frac{\sqrt{n}\big(\widehat G_n(Y_i) - G_0(Y_i)\big)}{G^2_0(Y_i)} \big(1+o_\PP(1)\big) \\
	& =  \frac{1}{G_0(Y_i)} \frac{1}{\sqrt{n}} \sum^n_{j=1} \int^B_0 \e1_{Y_i \geq s} \frac{dM_j^{\cal C}(s)}{y(s)} \big(1+o_\PP(1)\big).
\end{split}
\end{equation}
Note that for the martingale representation of  $(G_0 - \widehat G_n)/G_0$, the reader can see \cite{Fleming.Harrington.91}, page 97 or \cite{Chen.Jewell.05}. 
Thus
\begin{align*}
	{\cal Q}_{3n}-{\cal Q}_{2n} &= \sum^n_{i=1} \delta_i \bigg(  \frac{1}{\widehat G_n(Y_i)} - \frac{1}{G_0(Y_i)}\bigg) \big( \rho_\tau(\varepsilon_i - n^{1/2} \eX^\top_i \eu) - \rho_\tau (\varepsilon_i)\big) \\
	& =\frac{1}{n}  \sum^n_{i=1} \frac{\delta_i}{G_0(Y_i)}  \bigg(g_\tau (\varepsilon_i ) \frac{\eX^\top_i \eu}{ \sqrt{n}} +\frac{h_\tau(\varepsilon_i)}{2} \eu \frac{\eX_i \eX^\top_i}{n} \eu \bigg) \sum^n_{j=1} \int^B_0 \e1_{Y_i \geq s} \frac{d M^{\cal C}_j(s)}{y(s)}  \big(1+o_\PP(1)\big) \\
	& =\frac{1}{n}  \sum^n_{i=1} \frac{\delta_i}{G_0(Y_i)}g_\tau (\varepsilon_i ) \frac{\eX^\top_i \eu}{ \sqrt{n}} \sum^n_{j=1} \int^B_0 \e1_{Y_i \geq s}\frac{d M^{\cal C}_j(s)}{y(s)}  \big(1+o_\PP(1)\big).
\end{align*}
For the second equality we used relation (\ref{rell}), together with the fact that  $\| \eX \| < \infty$ of assumption (A5) and with  relation (\ref{GuZou}) for $t=n^{-1/2} \eX_i^\top \eu$. 
Thus, by the law of large numbers, we have:
\begin{equation}
	\label{Q32n}
	{\cal Q}_{3n}-{\cal Q}_{2n} =\frac{\eu^\top}{\sqrt{n}}\sum^n_{j=1}\int^B_0\frac{\ek(s)}{y(s)}dM^{\cal C}_j(s)\big(1+o_\PP(1)\big).
\end{equation}
By the martingale CLT we have:
\begin{equation}
	\label{Q32nc}
	\frac{\eu^\top}{ \sqrt{n}}\sum^n_{j=1}\int^B_0\frac{\ek(s)}{y(s)}dM^{\cal C}_j(s)   \overset{\cal L} {\underset{n \rightarrow \infty}{\longrightarrow}} \eu^\top \textbf{W}_2, 
\end{equation}
with the random $p$-vector  $\textbf{W}_2 \sim {\cal N}_p(\oo_p, \textbf{S}_2)$. \\
Hence, we have shown, by relations (\ref{QQQ}), (\ref{Q1n}), (\ref{Q32n}) and (\ref{Q32nc}), that:
\[
Q_n(\widehat G_n, \ebo+n^{-1/2} \eu) - Q_n(\widehat G_n, \ebo)=\bigg(\eu^\top \textbf{W}_1  +\frac{\eE_\varepsilon \big[h_\tau(\varepsilon)\big]}{2} \eu^\top \eE_\eX[\eX \eX^\top] \eu+\eu^\top \textbf{W}_2\bigg)\big(1+o_\PP(1)\big),
\]
where, for $c=\|\eu \|$ large enough, we have that $ \eE_\varepsilon \big[h_\tau(\varepsilon)\big]  \eu^\top \eE_\eX[\eX \eX^\top] \eu $ dominates $\eu^\top (\textbf{W}_1+\textbf{W}_2)$. This implies that for all $\epsilon >0$, for large $n$, we have:
\[
\PP \bigg[ \inf_{\| \eu \|=c} \big[ Q_n(\widehat G_n, \ebo+n^{-1/2} \eu)- Q_n(\widehat G_n, \ebo)\big]>0\bigg]  \geq 1-\epsilon,
\]
which involves that $\widetilde \eb_n -\ebo =O_\PP(n^{-1/2})$ and statement \textit{(i)} is proved.\\
\textit{(ii)} In view of statement \textit{(i)}, we have that the minimizer in $\eu$ of $Q_n(\widehat G_n, \ebo+n^{-1/2} \eu)$ is the solution of the following system of $p$ equations:
\[
\textbf{W}_1+\eE_\varepsilon \big[h_\tau(\varepsilon)\big] \eE_\eX [\eX \eX^\top]  \eu +\textbf{W}_2=\oo_p.
\]
The  solution of this system is  $\eu=-(\textbf{W}_1+\textbf{W}_2) \eE^{-1}_\varepsilon \big[h_\tau(\varepsilon)\big] \eE_\eX [\eX \eX^\top]^{-1}$.  On the other hand, since $\textbf{W}_1+\textbf{W}_2 \sim {\cal N}_p(\oo_p,\textbf{S}_1+\textbf{S}_2)$, then $\eu \sim {\cal N}_p \big(\oo_p, \textbf{S}^{-1}_3(\textbf{S}_1+\textbf{S}_2)\textbf{S}^{-1}_3\big)$ and statement \textit{(ii)} is proved.
\hspace*{\fill}$\blacksquare$  \\

\noindent {\bf Proof of  Theorem \ref{Theorem2.2 Tang}.}
For a survival function $G$ and a parameter $p$-vector $\eb$, let us consider the following random process:
\[
{\cal T}_n(G,\eb) \equiv Q_n(G,\eb) - Q_n(G,\ebo) +\lambda_n \sum^p_{j=1} \widehat{\omega}_{n,j} \big( |\beta_j|-|\beta^0_j|\big).
\]
We can write
\begin{equation}
	\label{dec1}
{\cal T}_n(\widehat G_n,\eb) = {\cal R}_n(\widehat G_n,\eb)  +\lambda_n \sum^p_{j=1} \widehat{\omega}_{n,j} \big( |\beta_j|-|\beta^0_j|\big),
\end{equation}
with ${\cal R}_n(\widehat G_n,\eb)$ defined in the proof of Theorem  \ref{Theorem2.1 Tang}. Consider the $p$-vector $\eu=(u_1, \cdots , u_p)$. For the second term of the right-hand side of relation (\ref{dec1}) we have with probability one  that  
\begin{align*}
	\sum^p_{j=1} \widehat{\omega}_{n,j} \big( |\beta^0_j+n^{-1/2} u_j|-|\beta^0_j|\big)
&  	=\sum^q_{j=1}\widehat{\omega}_{n,j} \big( |\beta^0_j+n^{-1/2} u_j|-|\beta^0_j|\big) + \sum^p_{j=q+1}\widehat{\omega}_{n,j} n^{-1/2}|u_j| \\
& \geq \sum^q_{j=1}\widehat{\omega}_{n,j} \big( |\beta^0_j+n^{-1/2} u_j|-|\beta^0_j|\big) .
\end{align*}
 On the other hand, using Theorem \ref{Theorem2.1 Tang} we have that $\widehat{\omega}_{n,j} =O_\PP(1)$ for all $j \in \{1 , \cdots, q\}$,  which implies with probability converging to one that:
\[
	\sum^p_{j=1} \widehat{\omega}_{n,j} \big( |\beta^0_j+n^{-1/2} u_j|-|\beta^0_j|\big) \geq  - c \sum^q_{j=1} n^{-1/2}| u_j|. 
\] 
Then, since $n^{-1/2} \lambda_n=O_\PP(1)$,  we have:
\begin{equation}
	\label{eq1b}
	\lambda_n  \sum^p_{j=1} \widehat{\omega}_{n,j} \big( |\beta^0_j+n^{-1/2} u_j|-|\beta^0_j|\big) \geq -O_\PP(\|\eu \|).
\end{equation}
On the other hand, by the proof of Theorem \ref{Theorem2.1 Tang}, for $\|\eu \|$ large enough, we have for the first term of the right-hand side of relation (\ref{dec1}):
\begin{equation}
	\label{eq2b}
	Q_n(\widehat G_n,\ebo+n^{-1/2} \eu) - Q_n(\widehat G_n,\ebo)=\frac{\eE_\varepsilon[h_\tau(\varepsilon)]}{2} \eu^\top \eE_\eX[\eX \eX^\top] \eu \big(1+o_\PP(1)\big).
\end{equation}
From relations (\ref{dec1}), (\ref{eq1b}) and (\ref{eq2b}) we obtain 
\[
{\cal T}_n\big(\widehat G_n,\ebo+n^{-1/2} \eu\big) =O_\PP(\|\eu \|^2),
\]
for $\|\eu \|$ large enough. Thus, $\widehat\eb_n -\ebo=O_\PP(n^{-1/2})$, which implies the theorem.
\hspace*{\fill}$\blacksquare$  \\

\noindent {\bf Proof of  Theorem \ref{Theorem2.3 Tang}}
\textit{(i)} Taking into account the convergence rate of $\widehat \eb_n$ towards $\ebo$ obtained by Theorem \ref{Theorem2.2 Tang}, let us consider the following sets of parameters ${\cal V}(\ebo) \equiv \big\{ \eb ; \| \eb-\ebo \| \leq c n^{-1/2} \big\}$ and ${\cal W}_n \equiv \big\{\eb \in {\cal V}(\ebo) ;  \| \eb_{{\cal A}^c} \| >0  \big\}$. Theorem \ref{Theorem2.2 Tang} implies that  $\widehat\eb_n$ belongs to ${\cal V}(\ebo)$ with probability converging to 1 when $n \rightarrow \infty$. In order to prove the theorem, we will first show that $\PP[\widehat \eb_n \in {\cal  W}_n]{\underset{n \rightarrow \infty}{\longrightarrow}}0$. \\
Recall that the true parameter is  $\ebo=\big(\eb^0_{\cal A}, \oo_{|{\cal A}^c|}\big)$. Then, in order to show the sparsity of $\widehat \eb_n$,  we consider two parameters $\eb^{(1)} \equiv \big( \eb^{(1)}_{\cal A}, \eb^{(1)}_{{\cal A}^c}\big) \in {\cal V}(\ebo) \setminus {\cal W}_n $ , $\eb^{(2)} \equiv \big( \eb^{(2)}_{\cal A}, \eb^{(2)}_{{\cal A}^c}\big) \in {\cal W}_n$, with $\eb^{(1)}_{{\cal A}^c} =\oo_{|{\cal A}^c|}$ and $ \eb^{(2)}_{\cal A} =  \eb^{(1)}_{\cal A}$.
Let's calculate 
\begin{equation}
	\label{TL}
	{\cal T}_n(\widehat G_n,\eb^{(2)}) - {\cal T}_n(\widehat G_n,\eb^{(1)}) = Q_n(\widehat G_n,\eb^{(2)}) - Q_n(\widehat G_n,\eb^{(1)})+\lambda_n \sum^p_{j=1} \widehat{\omega}_{n,j} \big( |\beta^{(2)}_j| - |\beta^{(1)}_j|\big).
\end{equation}
Then, for the first term of the right-hand side of relation (\ref{TL}),  taking in relation (\ref{GuZou}), $e=\varepsilon_i$ and $t=\eX^\top_i(\eb^{(2)}- \ebo)$ or $t=\eX^\top_i(\eb^{(1)}- \ebo)$ and using assumption (A5), we obtain:
\begin{align*}
	&Q_n(\widehat G_n,\eb^{(2)}) - Q_n(\widehat G_n,\eb^{(1)}) =\sum^n_{i=1} \frac{\delta_i}{\widehat G_n(Y_i)} \bigg( \rho_\tau \big( \log(Y_i) -\eX_i^\top \eb^{(2)} \big) - \rho_\tau \big(\log(Y_i) -\eX_i^\top \eb^{(1)} \big) \bigg)
	 \\
	&=\sum^n_{i=1} \frac{\delta_i}{\widehat G_n(Y_i)} \bigg\{\bigg( g_\tau (\varepsilon_i) \eX^\top_{i} \big(\eb^{(2)}_{{\cal A}} -\ebo_{\cal A},\eb^{(2)}_{{\cal A}^c} \big) +\frac{h_\tau(\varepsilon_i)}{2} \big( \eX^\top_{i} \big(\eb^{(2)}_{{\cal A}} -\ebo_{\cal A},\eb^{(2)}_{{\cal A}^c} \big)\big)^2 
	\\
	&  \qquad \qquad  +o_\PP\big( \eX^\top_{i} \big(\eb^{(2)}_{{\cal A}} -\ebo_{\cal A},\eb^{(2)}_{{\cal A}^c} \big)\big)^2\bigg)  \\
	&    - \bigg( g_\tau (\varepsilon_i) \eX^\top_{{\cal A},i} \big(\eb^{(2)} -\ebo \big)_{\cal A} +\frac{h_\tau(\varepsilon_i)}{2} \big( \eX^\top_{i,{\cal A}} \big(\eb^{(2)} -\ebo \big)_{\cal A}\big)^2 +o_\PP\big( \eX^\top_{{\cal A},i} \big(\eb^{(2)} -\ebo \big)_{\cal A}\big)^2\bigg) \bigg\} 
	\\
	&=\sum^n_{i=1} \frac{\delta_i}{G_0(Y_i)} \bigg\{\bigg(   g_\tau (\varepsilon_i) \eX^\top_{i} \big(\eb^{(2)}_{{\cal A}} -\ebo_{\cal A},\eb^{(2)}_{{\cal A}^c} \big) +\frac{h_\tau(\varepsilon_i)}{2} \big( \eX^\top_{i} \big(\eb^{(2)}_{{\cal A}} -\ebo_{\cal A},\eb^{(2)}_{{\cal A}^c} \big)\big)^2
	 \\
	& \qquad \qquad  +o_\PP\big( \eX^\top_{i} \big(\eb^{(2)}_{{\cal A}} -\ebo_{\cal A},\eb^{(2)}_{{\cal A}^c} \big)\big)^2\bigg)  
	\\
	& - \bigg( g_\tau (\varepsilon_i) \eX^\top_{{\cal A},i} \big(\eb^{(2)} -\ebo \big)_{\cal A} +\frac{h_\tau(\varepsilon_i)}{2} \big( \eX^\top_{{\cal A},i} \big(\eb^{(2)} -\ebo \big)_{\cal A}\big)^2 +o_\PP\big( \eX^\top_{{\cal A},i} \big(\eb^{(2)} -\ebo \big)_{\cal A}\big)^2\bigg)\bigg\} 
	\\
	&+\sum^n_{i=1} \delta_i \bigg( \frac{1}{\widehat G_n(Y_i)} -\frac{1}{G_0(Y_i)}\bigg) \bigg\{\bigg(  g_\tau (\varepsilon_i) \eX^\top_{i} \big(\eb^{(2)}_{{\cal A}} -\ebo_{\cal A},\eb^{(2)}_{{\cal A}^c} \big) +\frac{h_\tau(\varepsilon_i)}{2} \big( \eX^\top_{i} \big(\eb^{(2)}_{{\cal A}} -\ebo_{\cal A},\eb^{(2)}_{{\cal A}^c} \big)\big)^2 \\
	& \qquad \qquad  +o_\PP\big( \eX^\top_{i} \big(\eb^{(2)}_{{\cal A}} -\ebo_{\cal A},\eb^{(2)}_{{\cal A}^c} \big)\big)^2 \bigg)   
\\
	& - \bigg( g_\tau (\varepsilon_i) \eX^\top_{{\cal A},i} \big(\eb^{(2)} -\ebo \big)_{\cal A} +\frac{h_\tau(\varepsilon_i)}{2} \big( \eX^\top_{{\cal A},i} \big(\eb^{(2)} -\ebo \big)_{\cal A}\big)^2 +o_\PP\big( \eX^\top_{{\cal A},i} \big(\eb^{(2)} -\ebo \big)_{\cal A}\big)^2\bigg) \bigg\}
\end{align*}
Taking into account relation (\ref{rell}), this is equal to 
\begin{align*}
	&=\sum^n_{i=1} \frac{\delta_i}{G_0(Y_i)} \bigg\{ \bigg( g_\tau (\varepsilon_i) \eX^\top_{{\cal A}^c,i} \eb^{(2)}_{{\cal A}^c} 
	+\frac{h_\tau(\varepsilon_i)}{2} \big( \eX^\top_{i} \big(\eb^{(2)}_{{\cal A}} -\ebo_{\cal A},\eb^{(2)}_{{\cal A}^c} \big)\big)^2  +  o_\PP\big( \eX^\top_{i} \big(\eb^{(2)}_{{\cal A}} -\ebo_{\cal A},\eb^{(2)}_{{\cal A}^c} \big)\big)  \bigg) 
\\
	& \qquad - \bigg( \frac{h_\tau(\varepsilon_i)}{2} \big( \eX^\top_{{\cal A},i} \big(\eb^{(2)} -\ebo \big)_{\cal A}\big)^2 +o_\PP\big( \eX^\top_{{\cal A},i} \big(\eb^{(2)} -\ebo \big)_{\cal A}\big)\bigg) \bigg\} 
\end{align*}
\begin{align*}
	&+\sum^n_{i=1}  \frac{\delta_i}{G_0(Y_i)} \bigg\{  \bigg( g_\tau (\varepsilon_i) \eX^\top_{i,{\cal A}^c} \eb^{(2)}_{{\cal A}^c} +
	\frac{h_\tau(\varepsilon_i)}{2} \big( \eX^\top_{i} \big(\eb^{(2)}_{{\cal A}} -\ebo_{\cal A},\eb^{(2)}_{{\cal A}^c} \big)\big)^2 +o_\PP\big( \eX^\top_{i} \big(\eb^{(2)}_{{\cal A}} -\ebo_{\cal A},\eb^{(2)}_{{\cal A}^c} \big)\big)^2 \bigg) 
	\\
	&- \bigg( \frac{h_\tau(\varepsilon_i)}{2} \big( \eX^\top_{{\cal A},i} \big(\eb^{(2)} -\ebo \big)_{\cal A}\big)^2 +o_\PP\big( \eX^\top_{{\cal A},i} \big(\eb^{(2)} -\ebo \big)_{\cal A}\big)  \bigg) \bigg\}  \frac{1}{n} \sum^n_{j=1}\int^B_0 \e1_{Y_i \geq s} \frac{dM_j^{\cal C}(s)}{y(s)}\big(1+o_\PP(1)\big)\\
	& \equiv \SSS_1 +\SSS_2 \equiv (S_{11}-S_{12}) + (S_{21} - S_{22}).
\end{align*}
We analyze in the following each terms $S_{11}$, $S_{12}$, $S_{21}$ and $S_{22}$. Let's start with $S_{11}$. 
Since, by the CLT 
\[
n^{-1/2} \sum^n_{i=1} \frac{\delta_i}{G_0(Y_i)}  g_\tau(\varepsilon_i) \eX_{{\cal A}^c,i}  \overset{\cal L} {\underset{n \rightarrow \infty}{\longrightarrow}} {\cal N}_{|{\cal A}^c|}(\oo_q, \textbf{S}_{1,{\cal A}^c}),
\] 
then, we obtain
\begin{equation}
	\label{eq3}
	\sum^n_{i=1} \frac{\delta_i}{G_0(Y_i)}  g_\tau(\varepsilon_i) \eX^\top_{{\cal A}^c,i} \eb^{(2)}_{{\cal A}^c}=O_\PP(1).
\end{equation}
On the other hand, taking into account assumption (A5), we get: 
\begin{equation}
	\label{eq4}
	\frac{1}{n}\sum^n_{i=1} \frac{\delta_i}{G_0(Y_i)}  \frac{h_\tau(\varepsilon_i)}{2} \big( \eX^\top_{i} \big(\eb^{(2)}_{{\cal A}} -\ebo_{\cal A},\eb^{(2)}_{{\cal A}^c}\big)\big)^2=O_\PP\big( \| \eb^{(2)} -\ebo \|^2\big)=O_\PP(n^{-1}).
\end{equation}
Relations (\ref{eq3}) and (\ref{eq4}) imply that $S_{11}=O_\PP(1)$. We show similarly that $S_{12}=O_\PP(1)$ and then we obtain $\SSS_1=O_\PP(1)$.\\
For term $\SSS_2$, we first study $S_{21}$. Taking into account the definition of the random vector $\ek(s)$, we obtain:
\begin{equation}
\label{lp}
\frac{1}{n}\sum^n_{i=1}  \frac{\delta_i}{G_0(Y_i)} g_\tau(\varepsilon_i) \eX^\top_{{\cal A}^c,i}  \e1_{Y_i \geq s}=O_\PP(1) \overset{\PP} {\underset{n \rightarrow \infty}{\longrightarrow}} \ek_{{\cal A}^c}(s).
\end{equation}
By the martingale CLT we have
\[
\frac{1}{\sqrt{n}} \sum^n_{j=1} \int^B_0 \frac{ \ek_{{\cal A}^c}(s)}{y(s)} dM^{\cal C}_j(s) \overset{\cal L} {\underset{n \rightarrow \infty}{\longrightarrow}} \textbf{W}_{2,{\cal A}^c} \sim  {\cal N}_{|{\cal A}^c|}\big( \oo_{|{\cal A}^c|}, ...\big),
\]
which implies,  using Slutsky's Lemma, that, 
$$
n^{-1} \sum^n_{j=1} \int^B_0  \frac{\ek_{{\cal A}^c}(s)}{y(s)} dM^{\cal C}_j(s) \eb^{(2)}_{{\cal A}^c} =o_\PP(1).
$$
 Thus, using also relation (\ref{lp}), we obtain:
\[
\sum^n_{i=1}  \frac{\delta_i}{G_0(Y_i)} g_\tau(\varepsilon_i) \eX^\top_{{\cal A}^c,i} \eb^{(2)}_{{\cal A}^c} \frac{1}{n} \sum^n_{j=1}\int^B_0 \e1_{Y_i \geq s} \frac{dM_j^{\cal C}(s)}{y(s)} =O_\PP(1).
\]
Taking into account relation (\ref{eq4}) we obtain that $S_{21}=O_\PP(1)$. We show similarly that   $S_{22}=O_\PP(1)$. We have then $\SSS_2=O_\PP(1)$.\\
In conclusion, we showed that
\begin{equation}
	\label{eq6}
	Q_n(\widehat G_n,\eb^{(2)}) - Q_n(\widehat G_n,\eb^{(1)}) =O_\PP(1).
\end{equation}
On the other hand, taking into account Theorem \ref{Theorem2.1 Tang}\textit{(i)} and the assumption $n^{(\gamma -1)/2}\lambda_n  {\underset{n \rightarrow \infty}{\longrightarrow}} \infty$, we obtain:
\begin{equation}
	\label{eq7}
	\begin{split}
\lambda_n \sum^p_{j=1} \widehat{\omega}_{n,j} \big( |\beta^{(2)}_j| - |\beta^{(1)}_j|\big) & = \lambda_n \sum^p_{j=q+1} \widehat{\omega}_{n,j} |\beta^{(2)}_j| \\
 & 	\geq O_\PP \big( \lambda_n n^{\gamma/2} n^{-1/2}\big)  =O_\PP(\lambda_n n^{(\gamma -1)/2}) {\underset{n \rightarrow \infty}{\longrightarrow}} \infty.
	\end{split}
\end{equation}
Relations (\ref{eq6}) and (\ref{eq7}) imply that the penalty dominates in  ${\cal T}_n(\widehat G_n,\eb^{(2)}) - {\cal T}_n(\widehat G_n,\eb^{(1)})$ of relation (\ref{TL}), penalty which is the order  $\lambda_n n^{(\gamma -1)/2}{\underset{n \rightarrow \infty}{\longrightarrow}} \infty$. On the other hand, by a calculation similar to that of relation (\ref{eq6}), we obtain that:  ${\cal T}_n(\widehat G_n,\eb^{0}) -{\cal T}_n(\widehat G_n,\eb^{(1)})=	Q_n(\widehat G_n,\eb^0) - Q_n(\widehat G_n,\eb^{(1)})  =O_\PP(1)$. These imply that $\eb^{(2)}$ cannot be a minimum point of ${\cal T}_n(\widehat G_n,\eb )$. Thus, we have that $\PP[\widehat \eb_n \in {\cal  W}_n]{\underset{n \rightarrow \infty}{\longrightarrow}}0$ which implies $\PP[\widehat \eb_n \in {\cal A}^c \cap \widehat {\cal A}_n] {\underset{n \rightarrow \infty}{\longrightarrow}}0$, from where $\PP[  \widehat {\cal A}_n \subseteq {\cal A}] {\underset{n \rightarrow \infty}{\longrightarrow}}1$. Since $\widehat \eb_{n,{\cal A}}$ is consistent, then $\widehat \eb_{n,{\cal A}} \overset{\PP} {\underset{n \rightarrow \infty}{\longrightarrow}}\eb^0_{\cal A} \neq \oo_{|{\cal A}|}$, which implies $\PP[{\cal A} \subseteq \widehat {\cal A}_n ] {\underset{n \rightarrow \infty}{\longrightarrow}}1$. We have shown that $\PP[{\cal A} = \widehat {\cal A}_n ] {\underset{n \rightarrow \infty}{\longrightarrow}}1$, that is statement \textit{(i)}.\\
\textit{(ii)} In virtue of statement \textit{(i)}, we consider the parameter $\eb=(\eb_{\cal A},\eb_{{\cal A}^c})$ such that $\eb_{\cal A}=\eb^0_{\cal A}+n^{-1/2}\eu_{\cal A}$ and $\eb_{{\cal A}^c}=\oo_{|{\cal A}^c|}$, $\eu=(\eu_{\cal A}, \oo_{|{\cal A}^c|})$, with $\| \eu_{\cal A}\| \leq c \leq \infty$. Then, combining relation (\ref{GuZou})  with assumption (A5), we get:
\begin{align*}
	&{\cal T}_n \big( \widehat G_n, n^{-1/2}(\eu_{\cal A}, \oo_{|{\cal A}^c|}) \big)=Q_n\big(\widehat G_n, \ebo+n^{-1/2} \eu \big) - Q_n(\widehat G_n, \ebo)+\lambda_n \sum^q_{j=1} \big(|\beta^0_j +n^{-1/2}u_j|- |\beta^0_j|\big)\\
	& = \sum^n_{i=1} \frac{\delta_i}{\widehat G_n(Y_i)} \big( \rho_\tau(\varepsilon_i -n^{-1/2} \eX^\top_{{\cal A},i} \eu_{\cal A}) -\rho_\tau(\varepsilon_i)\big)+\lambda_n\sum^q_{j=1} \widehat{\omega}_{n,j} \frac{\sign(\beta^0_j)u_j}{\sqrt{n}}\\
	&= \bigg\{ \sum^n_{i=1} \frac{\delta_i}{  G_0(Y_i)} \bigg( g_\tau(\varepsilon_i) \frac{\eX^\top_{{\cal A},i} \eu_{\cal A}}{ \sqrt{n}} +\frac{h_\tau(\varepsilon_i)}{2}\bigg(\frac{\eX^\top_{{\cal A},i} \eu_{\cal A}}{ \sqrt{n}}\bigg)^2+o_\PP\bigg(\frac{1}{n}\bigg)\bigg)\\
	& \qquad  +\frac{1}{n} \sum^n_{i=1} \frac{\delta_i}{  G_0(Y_i)} g_\tau(\varepsilon_i) \frac{\eX^\top_{{\cal A},i} \eu_{\cal A}}{ \sqrt{n}}\sum^n_{j=1}\int^B_0 \frac{\e1_{Y_i \geq s}}{y(s)} dM^{\cal C}_j(s) \bigg\}\big(1+o_\PP(1)\big) +\lambda_n \sum^q_{j=1} \widehat{\omega}_{n,j} \frac{\sign(\beta^0_j)u_j}{\sqrt{n}},
\end{align*}
the $n^{-1}$ from the last line comes from Taylor's expansion of $ {1}/{\widehat G_n(Y_i)} -  {1}/{G_0(Y_i)}$ given by relation (\ref{rell}).\\
We recall that $\ek_{\cal A}(s) \equiv \lim_{n \rightarrow \infty} n^{-1}\sum^n_{i=1} {\delta_i}/{  G_0(Y_i)} g_\tau(\varepsilon_i) \eX^\top_{{\cal A},i} \e1_{Y_i \geq s}$. Then, by an approach similar to that used in the proof of Theorem \ref{Theorem2.1 Tang}\textit{(i)}, taking into account the fact that  $\eE_{\cal C}[\delta_i]=G_0(Y_i)$, we have:
\begin{align*}
	& {\cal T}_n(\widehat G_n, \ebo+\frac{\eu}{\sqrt{n}})=\frac{1}{\sqrt{n}} \sum^n_{i=1} \frac{\delta_i}{  G_0(Y_i)} g_\tau(\varepsilon_i)  \eX^\top_{{\cal A},i} \eu_{\cal A}+ \frac{\eE_\varepsilon[h_\tau(\varepsilon)]}{2}\eu^\top_{\cal A} \eE_\eX[\eX_{{\cal A}}\eX^\top_{{\cal A}}]\eu_{\cal A}\\
	& \qquad  +\frac{\eu^\top_{\cal A}}{\sqrt{n}} \sum^n_{j=1}\int^B_0 \frac{\ek_{\cal A}(s)}{y(s)}dM^{\cal C}_j(s)+\lambda_n \sum^q_{j=1} \widehat{\omega}_{n,j} \frac{\sign(\beta^0_j)u_j}{\sqrt{n}}+o_\PP(1)
\end{align*}
\begin{align*}
	& =\frac{\eE_\varepsilon[h_\tau(\varepsilon)]}{2}\eu^\top_{\cal A} \eE_\eX[\eX_{{\cal A}}\eX^\top_{{\cal A}}]\eu_{\cal A} +\frac{\eu^\top_{\cal A}}{\sqrt{n}} \sum^n_{i=1}\bigg(\frac{\delta_i}{  G_0(Y_i)} g_\tau(\varepsilon_i)  \eX_{{\cal A},i}+\int^B_0 \frac{\ek_{\cal A}(s)}{y(s)}dM^{\cal C}_i(s) \bigg)\\
	& \qquad +\lambda_n \sum^q_{j=1} \widehat{\omega}_{n,j} \frac{\sign(\beta^0_j)u_j}{\sqrt{n}}+o_\PP(1)\\
	& =\frac{\eE_\varepsilon[h_\tau(\varepsilon)]}{2}\eu^\top_{\cal A} \eE_\eX[\eX_{{\cal A}}\eX^\top_{{\cal A}}]\eu_{\cal A} + \eu^\top_{\cal A} \textbf{W}_3+l_0 \widehat \eo_{n,{\cal A}} \sign(\eb^0_{\cal A})^\top \eu_{\cal A},
\end{align*}
with the random $|{\cal A}|$-vector $\textbf{W}_3 \sim {\cal N}_{|{\cal A}|}(\oo_{|{\cal A}|}, \textbf{V}_3)$ and the deterministic $|{\cal A}|$-vector $\sign(\eb^0_{\cal A})\equiv \sign(\beta^0_j)_{j \in {\cal A}}$.\\
Then, the minimizer in $\eu$ of ${\cal T}_n(\widehat G_n, \ebo+n^{-1/2} \eu ) $ is the solution of the following system of $|{\cal A}|$ equations: 
$$  
\frac{\partial {\cal T}_n(\widehat G_n, \ebo+n^{-1/2} \eu )}{\partial \eu_{\cal A}}=\oo_{|{\cal A}|},
$$
that is:
\[
\eE_\varepsilon[h_\tau(\varepsilon)] \eE_\eX[\eX_{{\cal A}}\eX^\top_{{\cal A}}]\eu_{\cal A} + \textbf{W}_3+l_0 \widehat \eo_{n,{\cal A}} \sign(\eb^0_{\cal A}) =\oo_{|{\cal A}|},
\]
which implies that the solution is:
\[
\eu_{\cal A}  = - \eE^{-1}_\varepsilon[h_\tau(\varepsilon)]\eE_\eX[\eX_{{\cal A}}\eX^\top_{{\cal A}}]^{-1} \bigg(\textbf{W}_3+l_0 \widehat \eo_{n,{\cal A}} \sign(\eb^0_{\cal A}) \bigg).
\]
Then, taking into account the distribution of the random $|{\cal A}|$-vector $\textbf{W}_3$ we obtain:
\[
\widehat \eu_{n,{\cal A}}  \equiv \sqrt{n} \big( \widehat{\eb}_{n,{\cal A}} - \eb^0_{\cal A}\big) \overset{\cal L} {\underset{n \rightarrow \infty}{\longrightarrow}}  {\cal N}_{|{\cal A}|} \bigg(- \eE^{-1}_\varepsilon[h_\tau(\varepsilon)] l_0  {\eo^0_{\cal A}}^\top \sign(\eb^0_{\cal A})\eE_\eX[\eX_{{\cal A}}\eX^\top_{{\cal A}}]^{-1}, \eXL \bigg) 
\]
and the proof of statement \textit{(ii)} is finished. 
\hspace*{\fill}$\blacksquare$

\end{document}